\documentclass[11pt,amsfonts,fullpage]{amsart}

\usepackage{fullpage}

\newtheorem{theorem}{Theorem}[section]

\newtheorem{corollary}[theorem]{Corollary}

\newtheorem{proposition}[theorem]{Proposition}

\newtheorem{remark}[theorem]{Remark}

\def\endproof{\qed \medskip}
\def\blacksquare{\hbox to .60em{\vrule width .60em height .60em}}

\begin{document}

\title[ ]{A Survey of Einstein metrics on 4-manifolds}

\author[ ]{Michael T. Anderson}

\thanks{
Partially supported by NSF Grant DMS 0604735}

\maketitle

\abstract
We survey recent results and current issues on the existence and uniqueness 
of Einstein metrics on 4-manifolds. A number of open problems and conjectures 
are presented during the course of the discussion. 
\endabstract

\setcounter{section}{0}

\section{Introduction.}
\setcounter{equation}{0}

 The Einstein equations
\begin{equation}\label{e1.1}
Ric_{g} = \lambda g, \ \ \lambda \in {\mathbb R},
\end{equation}
for a Riemannian metric $g$ are the simplest and most natural set of 
equations for a metric on a given compact manifold $M$. Historically 
these equations arose in the context of Einstein's general theory of 
relativity, where the metric $g$ is of Lorentzian signature. However, 
over the past several decades there has also been much mathematical interest 
in Einstein metrics of Riemannian signature on compact manifolds, especially 
in low dimensions, and in particular in relation to the topology of the 
underlying manifold. 

 A strong motivation for this comes from the understanding developed 
in dimension 2 and more recently in dimension 3. To explain this, there is 
a complete classification of compact oriented 2-manifolds by the Euler 
characteristic, originally obtained by purely topological methods through 
work of M\"obius, Dehn, Heegard and Rado. This classification was 
later reproved via the Poincar\'e-Koebe uniformization theorem for surfaces, 
i.e.~any compact oriented surface carries a metric of constant curvature. The 
structure of such metrics then easily gives the full list of possible 
topological types of such surfaces, (and much more). 

 It has long been a goal of mathematicians to prove a similar classification 
of compact oriented 3-manifolds. Thurston \cite{Th2} realized that the key to 
this should be in studying the possible geometric stuctures on 
3-manifolds; the most important such structures are again the constant 
curvature, i.e.~Einstein, metrics. However, in contrast to surfaces, 
most 3-manifolds do not admit an Einstein metric \eqref{e1.1}; instead one 
had a simple, and conjecturally complete, list of well understood 
obstructions. A general 3-manifold should decompose into a collection 
of domains, each of which carries a natural geometry, the most important 
geometry being that of Einstein metrics. 

 The recent solution of the Thurston geometrization conjecture by 
Perelman \cite{Pe1}-\cite{Pe3} and Hamilton \cite{Ha1}-\cite{Ha2} has 
accomplished this goal of completely classifying all 3-manifolds. This 
has been obtained by understanding how to obtain solutions to the 
Einstein equations via the parabolic analogue of \eqref{e1.1}, namely 
the Ricci flow. 

 Ideally, one would like to carry out a similar program in dimension 4. 
However, (as with the passage from 2 to 3 dimensions), the world of 
4-manifolds is much more complicated than lower dimensions. One encounters 
a vast variety of exotic smooth structures, there are severe complications 
in understanding the fundamental group, and so on. Moreover, there are no 
canonical local models of Einstein metrics, and even when such metrics 
exist, the tie of the global geometry of these metrics with the underlying 
topology remains currently poorly understood. In fact, as discussed eloquently 
by Gromov in \cite{Gr3}, the ''dream'' of trying to understand $4$-manifolds 
via Einstein or other canonical metrics may well be impossible to realize. On 
the other hand, one should keep in mind that the deepest understanding to date 
of smooth $4$-manifolds comes via the geometry of connections or gauge fields, 
in the theories developed by Yang-Mills, Donaldson and Seiberg-Witten. How far 
such theories can be carried over to metrics, (the gravitational field), 
remains to be seen. 

\medskip

  This paper is an introductory survey of basic results to date on the existence, 
uniqueness, and structure of moduli spaces of Einstein metrics on 4-manifolds. 
Many interesting topics have been omitted or presented only briefly, due partly to 
the limits of the author's knowledge and taste, and partly to keep the article at 
a reasonable length. For this reason, we have excluded all discussion of Einstein 
metrics in higher dimensions. Also, any area of research is only as vital as the 
interest of significant open questions and problems. Accordingly, we present a number 
of such open problems throughout the paper, some of which are well-known and others 
less so.

\section{Brief Review: 4-manifolds, complex surfaces and Einstein metrics.}
\setcounter{equation}{0}

  In order to set the stage for the more detailed discussion to follow, 
in this section we present a very brief overview of the topology of 4-manifolds 
and the classification of complex surfaces. 

\medskip

{\bf I.} {\it Topology of 4-manifolds}.

  Throughout the paper, $M$ will denote a compact, oriented 4-manifold unless 
otherwise indicated. Together with the fundamental group $\pi_{1}(M)$, the 
most important topological datum of a 4-manifold is the cup product or 
intersection pairing,
\begin{equation}\label{e2.1}
I: H_{2}(M, {\mathbb Z})\otimes H_{2}(M, {\mathbb Z}) \rightarrow  
{\mathbb Z}.
\end{equation}
By Poincar\'e duality, $I$ is a symmetric, non-degenerate bilinear form when 
$M$ is simply connected, (so that $H_{2}(M, {\mathbb Z})$ is torsion-free). 
The algebraic classification of such pairings $(A, I)$ over ${\mathbb Z}$, where 
$A$ is an abelian group, starts with a few simple invariants. Thus, the rank 
$rank(I)$ and index $\tau(I)$ are defined by $rank(I) = p(I) + n(I)$
and $\tau(I)  = p(I) - n(I)$, where $p$ and $n$ are the maximal dimensions 
on which the form $I\otimes {\mathbb R}$ is positive or negative definite. 
The form $I$ is called even if $I(a,a) \equiv 0$ mod 2, for all $a\in A$; 
otherwise $I$ is odd. The algebraic classification then states that 
any indefinite pairing $(A, I)$ is determined up to isomorphism by its rank, 
index and parity, and has a simple standard normal form. If $(A, I)$ is 
definite, there is a huge number of possibilities; for example, there are 
more than $10^{51}$ inequivalent definite forms of rank $40$. There 
are only finitely many definite pairings of a given rank although a complete 
classification remains to be determined. 

 One has the obvious relations $|\tau(I)| \leq rank(I)$, $\tau(I) \equiv  
rank(I)$ (mod 2). In addition, if $I$ is even, then $\tau(I) \equiv 0$ 
(mod 8). Modulo these relations, all values of rank, index and parity are 
possible algebraically. 

 Now for a given compact oriented 4-manifold $M$ and for $I$ as in 
\eqref{e2.1}, set $b_{2}^{+}(M) = p(I)$, $b_{2}^{-}(M) = n(I)$ and 
$\tau(M) = \tau (I)$; $\tau(M)$ is called the signature of $M$. A simply 
connected 4-manifold $M$ is spin if $I$ is even, non-spin otherwise. A first 
basic result, due to J.H.C. Whitehead, is that the homotopy type of a simply 
connected 4-manifold is completely determined by the intersection pairing 
\eqref{e2.1}; thus, if $M$ and $N$ are simply connected, then they are 
homotopy equivalent if and only if $I(M) \simeq I(N)$. 

 In a remarkable work, Freedman \cite{F} established the classification at the 
next level, i.e.~up to homeomorphism. Thus, every pairing $(A, I)$ occurs 
as the intersection form of a simply connected topological 4-manifold. 
If $M$ is spin, then the homeomorphism type is unique, while if $M$ is 
non-spin, there are exactly two homeomorphism types; one and only one 
of these is stably smoothable, in that the product with ${\mathbb R}$ has 
a smooth structure. 

  Passing next to the category of smooth manifolds, the Rochlin theorem gives 
the further restriction $\tau(M) \equiv 0$ (mod 16) when $M$ is spin. Shortly 
after Freedman's work, Donaldson \cite{D1} proved the amazing result that if 
the intersection pairing $I(M)$ of a smooth, simply connected 4-manifold is 
definite, then $I(M)$ is a diagonal form over ${\mathbb Z}$ and hence, by 
Freedman's result, $M$ is homeomorphic to a connected sum $\#k{\mathbb C}{\mathbb P}^{2}$. 
This led to the existence of distinct smooth structures on 4-manifolds, and together 
with the Seiberg-Witten equations, has led to a vast zoo of possible smooth structures 
on a given underlying topological 4-manifold. In fact, it is possible that every smoothable 
4-manifold always has infinitely many distinct smooth structures, cf.~\cite{FS} 
and references therein for further details. 

  While there have been some generalizations of the results of Freedman and 
Donaldson to 4-manifolds with non-trivial fundamental group, comparatively little 
is known when the fundamental group is unrestricted. This is partly due to  
well-known undecidability issues; for instance, there is no algorithm to 
determine whether a general compact 4-manifold is simply connected, or 
whether two 4-manifolds are homeomorphic. 

\medskip

{\bf II.} {\it Complex surfaces}.

 Throughout this subsection, $M$ will be a compact complex surface. 
We briefly review the Kodaira classification of surfaces. This derives 
from the structure of the canonical line bundle $K$ over $M$, i.e.~the 
bundle of holomorphic 2-forms on $M$. In local holomorphic coordinates 
$(z^{1},z^{2})$, sections of $K$ have the form $f(z_{1},z_{2})dz^{1}\wedge 
dz^{2}$ with $f$ holomorphic. Let $K^{n}$ denote the $n$-fold tensor product 
of $K$, a line bundle with local sections of the form $f(z_{1},z_{2})(dz^{1}
\wedge dz^{2})^{n}$. The $n^{th}$ plurigenus of $M$ is the dimension of the 
space of holomorphic sections of $K^{n}$, $P_{n}(M) = dim_{{\mathbb C}}H^{0}
(M, K^{n})$.

 If $P_{n}(M) = 0$ for all $n > 0$, $M$ is said to have Kodaira 
dimension $\kappa (M) = -\infty$. Otherwise, one has $P_{n}(M) = 
O(n^{a})$, for $a = 0$, $1$ or $2$, and the smallest such a defines 
$\kappa (M)$. Thus $\kappa (M) \in \{-\infty, 0, 1, 2\}$. (Roughly 
speaking, $\kappa(M)$ is the dimension of the image of $M$ under the 
Kodaira map). 

 Next, recall the process of blowing up and down. Given any $M$ and 
$p \in M$, the blow-up $\hat M$ of $M$ at $p$ is the complex surface 
obtained by replacing a ball $B \subset {\mathbb C}^{2}$ near $p$ by 
$\bar{\mathbb C}{\mathbb P}^{2}\setminus B$, where $\bar{\mathbb C}
{\mathbb P}^{2}$ is ${\mathbb C}{\mathbb P}^{2}$ with the opposite 
orientation. Thus $\hat M = M \#\bar{\mathbb C}{\mathbb P}^{2}$ 
topologically. There is a holomorphic map $\pi: \hat M \rightarrow M$ 
and a rational curve $E = {\mathbb C}{\mathbb P}^{1} \subset \hat M$ 
such that $\pi$ restricted to $\hat M\setminus E$ is a biholomorphism 
onto $M\setminus \{p\}$, so that $\pi$ contracts $E$ to $\{p\}$. One has 
$[E]\cdot  [E] = -1$, i.e. $E$ has self-intersection -1 in 
$\hat M$ and $K(\hat M) = \pi^{*}K(M) + E$, (as divisors on 
$\hat M$). It follows that $\kappa(M)$ is unchanged under blow ups. 

 Whenever a complex surface $M$ has an exceptional curve $E$, i.e.~a 
rational curve with self-intersection -1, this process may be inverted, 
i.e. $M$ may be blown down to remove the divisor $E$. A surface $M$ is 
minimal if $M$ has no such exceptional curves. Every surface may be 
blown down (not necessarily uniquely) to a minimal surface. The Kodaira 
classification then states the following: 

 $\bullet$ If $\kappa (M) = -\infty $ and $M$ is simply connected, then $M$ 
is rational; $M$ can be obtained by blowing up ${\mathbb C}{\mathbb P}^{2}$ a 
finite number of times and then blowing down a finite number of times. 
The underlying 4-manifold $M$ is diffeomorphic to ${\mathbb C}{\mathbb P}^{2}
\#k\bar{\mathbb C}{\mathbb P}^{2}$, $k \geq 0$, or ${\mathbb C}{\mathbb P}^{1}
\times {\mathbb C}{\mathbb P}^{1}$. If $M$ is not simply connected, then $M$ 
is the blow-up of a minimal ruled surface; there is a holomorphic map 
$\pi: M \rightarrow S$, where $S$ is a complex curve, with fibers 
${\mathbb C}{\mathbb P}^{1}$. 

 If $\kappa (M) \geq 0$, then $M$ can be uniquely blown down to a 
minimal surface $M_{\min}$, called the minimal model, and from now on, 
assume then that $M$ is minimal. 

 $\bullet$ If $\kappa (M) = 0$, then $M$ is a K3 surface or complex torus, 
or a finite quotient of one of these spaces. Any K3 surface is diffeomorphic 
to a quadric in ${\mathbb C}{\mathbb P}^{3}$. 

 $\bullet$ If $\kappa (M) = 1$, then $M$ is elliptic; there is a holomorphic 
map $\pi : M \rightarrow C$, where $C$ is a complex curve, with fibers a 
smooth curve of genus 1, (an elliptic curve), for almost all $p \in C$.

 $\bullet$ If $\kappa (M) = 2$, then $M$ is of general type. For example, 
hypersurfaces of degree $k \geq 5$ in ${\mathbb C}{\mathbb P}^{3}$ are 
of general type. 

\medskip

{\bf III.} {\it Einstein metrics}.

  When expressed in local coordinates, the Einstein equations \eqref{e1.1} 
form a complicated, quasi-linear system of PDE's for the metric $g = g_{ij}$. 
This system is not elliptic, due to the diffeomorphism invariance of the 
equations. However, for a suitable choice of local slice transverse to the action 
of the diffeomorphism group on the space of metrics, the restricted equations 
become elliptic. The simplest such choice locally is the harmonic coordinate 
gauge, where each of the coordinate functions $x^{i}$ is harmonic with respect 
to the given metric $g$, i.e.~$\Delta_{g}x^{i} = 0$. In such coordinates, 
the Einstein equations have the pleasant form
\begin{equation}\label{e2.2}
-{\tfrac{1}{2}}\Delta_{g}g_{ij} + Q_{ij}(g, \partial g) = \lambda g_{ij},
\end{equation}
where $Q$ is quadratic in the metric $g$ and its first derivatives. Thus, the 
system of Einstein equations locally can be viewed as a non-linear and coupled 
version of the equations for eigenfunctions of the Laplacian $\Delta$. The 
ellipticity of the system \eqref{e2.2} implies that Einstein metrics are 
$C^{\infty}$, (in fact real-analytic) in such harmonic coordinates. 

   In dimensions 2 and 3, Einstein metrics are of constant curvature, and so 
locally rigid; they are all locally isometric to domains in the space-forms 
of constant curvature. This is no longer the case in dimensions 4 and above. 
In fact, the space of local solutions, (i.e.~solutions defined on a ball), of 
the equations \eqref{e1.1} is infinite dimensional, (analogous to the scalar 
eigenfunction equation). 

  In dimension 4, Einstein metrics enjoy a variational characterization, for 
instance as metrics minimizing the $L^{2}$ norm of the curvature tensor $R$, 
(cf.~\eqref{e4.2} below), or as critical points of the volume-normalized 
total scalar curvature functional 
\begin{equation}\label{e2.3}
{\mathcal S}(g) = (vol_{g}(M))^{-1/2}\int_{M}s_{g}dV_{g},
\end{equation}
where $s_{g}$ is the scalar curvature of $g$. It would be very interesting if 
Einstein metrics could be constructed variationally, by using methods in the 
calculus of variations for either of the functionals above. For example, 
the solution of the Yamabe problem \cite{Sc,LP} gives the existence of a 
metric, called a Yamabe metric, minimizing ${\mathcal S}$ in the conformal 
class $[g]$ of any metric $g$ on $M$. One would then like to understand the 
existence of Einstein metrics on $M$ by understanding the limiting behavior 
of sequences of Yamabe metrics $g_{i}$, whose scalar curvature tends to its 
maximal value. Unfortunately, comparatively little is currently known about 
this approach.

\section{Constructions of Einstein metrics I.}
\setcounter{equation}{0}

 In this section we discuss the primary method of proving the existence of 
Einstein metrics on 4-manifolds, namely the existence of K\"ahler-Einstein 
metrics via the solution of the Calabi conjecture. We also list the 
remaining currently known Einstein metrics on 4-manifolds; see also 
\S 7 for further discussion. 

 To begin, suppose $(M, J)$ is a compact complex 4-manifold admitting a 
K\"ahler metric $g$, with K\"ahler form $\omega  = g(J, \cdot)$. For 
K\"ahler metrics, the Ricci form $\rho = Ric(J, )$ is given, (up to a 
factor of $i$), by the curvature form of the canonical line bundle $K$ 
over $M$. The Ricci form thus depends only on the complex structure $J$ 
and the volume form $\mu_{g}$ of the K\"ahler metric $g$ and
\begin{equation}\label{e3.1}
[\rho] = \frac{1}{2\pi}[c_{1}(M)] \in  H^{2}(M, {\mathbb C}).
\end{equation}
Consequently, if $g$ is in addition Einstein, satisfying \eqref{e1.1}, 
then one has
\begin{equation}\label{e3.2}
[c_{1}(M)] = 2\pi\lambda [\omega].
\end{equation}
In particular $c_{1}< 0$, $c_{1} = 0$ or $c_{1} > 0$ according to 
whether $\lambda < 0$, $\lambda = 0$ or $\lambda > 0$. This of 
course gives a strong restriction on the existence of K\"ahler-Einstein 
metrics; the first Chern class $c_{1}$ must be definite or identically 0, 
(in the sense that the associated symmetric bilinear form has these 
properties). Thus, ``most'' complex surfaces do not admit a 
K\"ahler-Einstein metric. 

 On the other hand, a basic reason why one is able to solve the existence 
problem for K\"ahler-Einstein metrics, as opposed to solving the existence 
problem for general Einstein metrics, is that via \eqref{e3.1}-\eqref{e3.2}, 
the existence of K\"ahler-Einstein metrics is tightly bound to the cohomology 
of $M$. In this way, the K\"ahler-Einstein equation can be reduced to a single 
scalar PDE, instead of the complicated full system of PDE's in, for 
instance, \eqref{e2.2}. 

 It is worth describing this remarkable reduction in more detail. A given 
K\"ahler metric $g$ has the form
$$g = \sum g_{k\bar j}dz^{k}d\bar z^{j}, $$
in local complex coordinates. The Ricci tensor is then given by
$$Ric_{k\bar l} = -\partial_{k}\partial_{\bar l}\log det(g_{m\bar{n}}), $$
and the equation $[\rho] = \lambda [\omega]$ implies there is a smooth 
function $F$ such that $Ric_{k\bar l} - \lambda g_{k\bar l} = 
\partial_{k}\partial_{\bar l}F$. To study the existence of a 
K\"ahler-Einstein metric with the same K\"ahler class as $\omega$, let 
$\widetilde g_{k\bar l} = g_{k\bar l} + \partial_{k}\partial_{\bar l}\phi$. 
Then the Einstein equation $Ric_{\widetilde g} = \lambda \widetilde g$ is 
equivalent to the equation
\begin{equation}\label{e3.3}
\frac{det (g_{k\bar l} + 
\phi_{k\bar l})}{det(g_{k\bar l})} = e^{-\lambda\phi + F}.
\end{equation}
To see this, taking $\partial_{k}\partial_{\bar l}$ of the logarithm of 
both sides of \eqref{e3.3} gives
$$-(\widetilde Ric_{k\bar l} - Ric_{k\bar l}) = 
-\lambda\partial_{k}\partial_{\bar l}\phi  + 
\partial_{k}\partial_{\bar l}F = -\lambda 
(\widetilde g_{k\bar l} - g_{k\bar l}) + (Ric_{k\bar l} - 
\lambda g_{k\bar l}) = -\lambda\widetilde g_{k\bar l} + 
Ric_{k\bar l}, $$
which implies $\widetilde Ric_{k\bar l} = 
\lambda\widetilde g_{k\bar l}$. Thus a K\"ahler-Einstein metric 
exists in the class $[\omega]$ if and only if there is a solution 
$\phi$ of the scalar complex Monge-Ampere equation \eqref{e3.3} with 
$g_{k\bar l} + \phi_{k\bar l}$ positive definite. 

 The existence of a solution of \eqref{e3.3} is proved by the method of 
continuity, solving \eqref{e3.3} on a curve $\phi_{t} = t\phi$ with a suitable 
choice of $F_{t}$. The initial set-up above shows that a solution exists for 
$t = 0$, and one proves that the set of $t\in [0,1]$ for which \eqref{e3.3} 
has a solution $\phi_{t}$ is both open and closed. The openness result is a 
straightforward consequence of the inverse function theorem. The main task 
is to prove closedness of the set of solutions; this requires rather difficult 
apriori estimates on the behavior of the solutions. 

 The basic results are as follows:

\begin{theorem}\label{t3.1}\cite{Au}, \cite{Y1}. A compact 
complex surface $(M, J)$ admits a K\"ahler-Einstein metric with $\lambda < 0$ if 
and only if $c_{1} < 0$. This occurs precisely for $(M, J)$ which are minimal 
surfaces of general type which contain no $(-2)$-curves, (rational curves of 
self-intersection $-2$). The K\"ahler-Einstein metric is uniquely determined 
by $(M, J)$. 
\end{theorem}

\begin{theorem}\label{t3.2} \cite{Y1}. A compact complex surface $(M, J)$ admits a 
K\"ahler-Einstein metric with $\lambda = 0$ if and only if $c_{1} = 0$. This 
occurs precisely when $(M, J)$ is finitely covered by a K3 surface or a complex 
torus. The K\"ahler-Einstein metric is uniquely determined by $(M, J)$ and a 
choice of K\"ahler class in the K\"ahler cone in $H^{1,1}(M, {\mathbb R})$.
\end{theorem}

\begin{theorem}\label{t3.3} \cite{Ti}. A compact complex surface $(M, J)$ admits 
a K\"ahler-Einstein metric with $\lambda > 0$ if and only if $c_{1}> 0$, and 
the Lie algebra of holomorphic vector fields is reductive. This occurs 
exactly on ${\mathbb C}{\mathbb P}^{2}, {\mathbb C}{\mathbb P}^{1}\times 
{\mathbb C}{\mathbb P}^{1}$ or the blow-up ${\mathbb C}{\mathbb P}^{2}
\#k\bar{\mathbb C}{\mathbb P}^{2}$, $3 \leq  k \leq 8$, of 
${\mathbb C}{\mathbb P}^{2}$ at $k$ points in general 
position. Again, the K\"ahler-Einstein metric is uniquely determined by 
$(M, J)$, a result of \cite{BM}. 
\end{theorem}

 These results give a complete understanding of K\"ahler-Einstein metrics on 
a given complex surface, at least concerning formal existence and uniqueness 
issues. Note that there are relatively few K\"ahler-Einstein metrics 
when $\lambda \geq 0$; most all have $\lambda < 0$, in analogy to 
the case of Riemann surfaces. 

\medskip

  The results above show that the moduli space ${\mathcal E}_{KE}$ of 
K\"ahler-Einstein metrics may be naturally identified with the moduli 
space ${\mathcal M}_{C}$ of complex structures on a given manifold $M$ 
when $c_{1}(J) < 0$ or $c_{1}(J) > 0$, and similarly with the space of 
K\"ahler classes over ${\mathcal M}_{C}$ when $c_{1}(J) = 0$. A 
fundamental issue of interest, raised in particular by Yau, cf.~\cite{Y2} 
for instance, is then to use this identification to study in more detail 
the structure of each of these moduli spaces. This is still basically an 
undeveloped area and much work remains to be done in this direction. 

 Another well-known problem is to understand these K\"ahler-Einstein 
metrics more explicitly; see for instance the recent work of Donaldson 
\cite{D2} and further references therein. 

  One would also like to answer simple uniqueness questions. For instance, 
are K\"ahler-Einstein metrics the unique Einstein metrics on a given 
(complex) 4-manifold $M$?

\medskip

 A further collection of Einstein metrics on 4-manifolds are the locally 
homogeneous metrics, given by left-invariant metrics on a compact locally 
homogeneous space $\Gamma \setminus G/H$. A complete list of these Einstein 
metrics, modulo finite covers, is:
$${\mathbb S}^{4}, \ \ {\mathbb R}^{4}/\Gamma , \ \ {\mathbb H}^{4}/\Gamma , 
\ \ {\mathbb C}{\mathbb P}^{2}, \ \ {\mathbb S}^{2}\times {\mathbb S}^{2}, 
\ \ {\mathbb C}{\mathbb H}^{2}/\Gamma , \ \ {\mathbb H}^{2}/\Gamma_{1}\times 
{\mathbb H}^{2}/\Gamma_{2}, $$
all with their canonical metrics. Observe that all such metrics are in 
fact locally symmetric. 

 Next, one may consider Einstein metrics which although not (locally) 
homogeneous, have a non-trivial or rather large isometry group. Thus, 
recall a Riemannian manifold $(M, g)$ is of cohomogeneity $k$ if the 
isometry group $(M, g)$ acts on $M$ with principal orbits of codimension 
$k$. By a well-known theorem of Bochner, any Einstein metric with a non-trival 
connected isometry group which does not split a Euclidean factor 
must have positive scalar curvature, $\lambda > 0$. 

 A very interesting and explicit Einstein metric of cohomogeneity 1 on 
${\mathbb C}{\mathbb P}^{2}\#\bar{\mathbb C}{\mathbb P}^{2}$ was found by Page 
\cite{Pa1}; this has the form
$$g = V^{-1}dr^{2} + V\sigma_{1}^{2} + f(\sigma_{2}^{2} + 
\sigma_{3}^{2}), $$
where $\{\sigma_{i}\}$ are the standard coframing of $S^{3} = SU(2)$ and 
$V = V(r)$, $f = f(r)$ are explicit functions. This metric has an isometric 
$U(2)$ action and, up to finite covers, is the only known cohomogeneity 1 Einstein 
metric on a compact 4-manifold; it is still an open problem whether there are 
any other cohomogeneity 1 Einstein metrics on compact 4-manifolds, cf.~\cite{Da}. 

 Recently, another very interesting Einstein metric was found by 
Chen-LeBrun-Weber \cite{CLW} on the 2-point blow-up of ${\mathbb C}{\mathbb P}^{2}$, 
i.e.~${\mathbb C}{\mathbb P}^{2}\#2\bar{\mathbb C}{\mathbb P}^{2}$. 
This metric is toric, and so has an isometric $S^{1}\times S^{1}$ action Both 
the Page metric and the Chen-LeBrun-Weber metric are 
Hermitian-Einstein, and conformal to K\"ahler metrics, but are not K\"ahler 
themselves. 

 Together with the results on K\"ahler-Einstein metrics above, we see 
that the rational surfaces ${\mathbb C}{\mathbb P}^{1}\times {\mathbb C}
{\mathbb P}^{1}$ and ${\mathbb C}{\mathbb P}^{2}\#k\bar{\mathbb C}
{\mathbb P}^{2}$, $0 \leq k \leq 8$, all admit Einstein metrics. 
In \S 4, we will see that ${\mathbb C}{\mathbb P}^{2}\#k\bar{\mathbb C}
{\mathbb P}^{2}$ does not admit an Einstein metric, for $k \geq 9$. 

 It would be interesting to understand the class of cohomogeneity 2 
Einstein metrics on compact 4-manifolds in general. In the case of 
static $S^{1}\times S^{1}$ actions, where the metric is a double warped 
product over a 2-dimensional base, a classical procedure due to Weyl, 
cf.~\cite{St}, gives a method of constructing local Ricci-flat metrics 
from an axisymmetric harmonic function on a domain in ${\mathbb R}^{3}$; 
these are the so-called static axisymmetric vacuum metrics in general 
relativity and are governed by a single linear (!) scalar equation. Can 
this procedure be generalized to situations where $\lambda > 0$ or 
$\lambda < 0$? Can one construct new global solutions of the Einstein 
equations by this or related procedures?; see \cite{St}, and \cite{CP} 
for instance for further discussion. 

  There has also been very little study of Einstein metrics having only 
isometric $S^{1}$ actions, although again these are important in general 
relativity; see \cite{Se} for some non-trivial results in this direction.

 A further construction of Einstein metrics, analogous to the Thurston 
theory of Dehn surgery on hyperbolic 3-manifolds will be discussed in 
\S 7. Together with the metrics listed above, this constitutes the 
{\it complete} list of known Einstein metrics. Clearly, there is much 
further territory to explore here. 

\section{Obstructions to Einstein metrics.}
\setcounter{equation}{0}

 In this section, we discuss a number of the known obstructions to the 
existence of Einstein metrics on 4-manifolds. In tandem with this, we 
also discuss several rigidity or uniqueness results for such metrics on 
a given manifold. 

 The most elementary obstruction comes from a simple observation of Berger. 

\begin{theorem} \label{t4.1} \cite{Be}. If $(M^{4}, g)$ is an Einstein metric, 
then
\begin{equation} \label{e4.1}
\chi (M) \geq  0, 
\end{equation}
with equality if and only if $(M^{4}, g)$ is flat. 
\end{theorem}

{\bf Proof:} The Chern-Gauss-Bonnet formula in dimension 4 reads
\begin{equation} \label{e4.2}
\chi (M) = \frac{1}{8\pi^{2}}\int_{M}|R|^{2} - |z|^{2} = 
\frac{1}{8\pi^{2}}\int_{M}|W|^{2} - {\tfrac{1}{2}}|z|^{2} + 
{\tfrac{1}{24}}s^{2},
\end{equation}
where $z = Ric - \frac{s}{4}g$ is the trace-free Ricci curvature $R$ is 
the full Riemann curvature and $W$ is the Weyl curvature. Einstein metrics 
are characterized by the condition that $z = 0$, and so the result follows 
immediately. 
{\endproof}

 In a certain sense then, at least one-half of 4-manifolds do not admit 
Einstein metrics. Simple examples include all circle bundles over 
3-manifolds $M^{3} \neq  T^{3}$, or products of surfaces of genus of 
non-equal sign, e.g.~$S^{2}\times \Sigma_{g}$, $g \geq  1$, or 
$T^{2}\times \Sigma_{g}$, $g \geq 2$, (or more generally such surface 
bundles over surfaces). This kind of argument, although of course very 
simple, is typical of many of the known obstructions to the existence of 
Einstein metrics. One finds (sharp) inequalities among characteristic 
numbers of $M$, for which equality implies a rigidity of Einstein metrics 
on the given space. Such rigidity then appears as the border between 
possible existence and non-existence. 

 For example, a strengthening of this argument was found by Hitchin and Thorpe 
(independently), by bringing in the Hirzebruch signature formula:
\begin{equation} \label{e4.3}
\tau (M) = \frac{1}{12\pi^{2}}\int_{M}|W_{+}|^{2} - |W_{-}|^{2}, 
\end{equation}
where $\tau(M)$ is the signature of $M$ and $W_{\pm}$ are the self-dual 
and anti-self-dual components of the Weyl tensor. The analysis of the 
equality case below is due to Hitchin. 

\begin{theorem} \label{t4.2}\cite{Hi1, Tj}. If $(M^{4}, g)$ is an Einstein 
metric, then
\begin{equation} \label{e4.4}
\chi (M) \geq  \frac{3}{2}|\tau (M)|, 
\end{equation}
with equality if and only if $(M^{4}, g)$ is flat, or $(M^{4}, g)$ is a 
Ricci-flat K\"ahler metric on the $K3$ surface (Calabi-Yau metric), or a 
quotient of such. 
\end{theorem}

{\bf  Sketch of Proof:} The inequality (4.4) is essentially an immediate 
consequence of the formulas (4.2) and (4.3). If equality holds, then $(M, g)$ 
has (anti)-self-dual curvature and zero scalar curvature, and so is Ricci-flat. 
If $M$ is simply connected, then an examination of the space $\Lambda^{+}(M)$ of 
self-dual 2-forms shows that $(M, g)$ is K\"ahler, so that $c_{1} = 0$, which 
gives the result. 
{\endproof}

 For example, for the rational surfaces $M_{k} = {\mathbb C}{\mathbb P}^{2}
\#k\bar{\mathbb C}{\mathbb P}^{2}$, one has $\chi(M_{k}) = 3+k$ and 
$\tau(M_{k}) = 1 - k$. Hence if $k \geq 9$, then $M_{k}$ does not admit 
an Einstein metric. As discussed in \S 3, $M_{k}$ does admit an Einstein metric 
for $k \leq 8$. 

  Note that the characteristic numbers that enter the results above are all 
homotopy invariant, and so the results apply to all possible smooth structures 
on a given topological 4-manifold. Gompf and Mrowka show in \cite{GM} that 
there are infinitely many distinct smooth structures within the homeomorphism 
class of the K3 surface; hence it follows immediately from Theorem 4.2 that 
none of these exotic smooth structures on K3 admits an Einstein metric. 

\medskip

  A strengthing of \eqref{e4.4} holds in the case of complex surfaces with 
$c_{1} < 0$. Namely, a standard result from complex surface theory gives 
$c_{2} = \chi$ and $c_{1}^{2} = 3\tau + 2\chi$, so that, for any compact 
complex surface,
$$\chi (M) -3\tau (M) = 3c_{2} - c_{1}^{2}.$$
Via Chern-Weil theory the terms $c_{2}$ and $c_{1}^{2}$ are given in terms of 
the curvature of Hermitian metrics on $M$ and by inspection, it easy to see 
that if $(M, J)$ admits a K\"ahler-Einstein metric, then
\begin{equation} \label{e4.5}
\chi (M) -3\tau (M) = 3c_{2} - c_{1}^{2} \geq  0, 
\end{equation}
with equality if and only if $(M, J)$ is biholomorphic to a complex 
hyperbolic space-form ${\mathbb C}{\mathbb H}^{2}/\Gamma$. By Theorem 3.1, 
the Bogomolov-Miyaoka-Yau inequality \eqref{e4.5} thus holds for any complex 
surface admitting a K\"ahler metric with $c_{1} < 0$. (This equation is 
uninteresting when $c_{1} = 0$ or $c_{1} > 0$). 

  Yau's proof above of the rigidity of ${\mathbb C}{\mathbb H}^{2}/\Gamma$ among 
complex surfaces and among K\"ahler-Einstein metrics does not extend to give 
a rigidity for (non-K\"ahler) Einstein metrics. In a series of papers, LeBrun 
has used the Seiberg-Witten theory to extend many of the rigidity and non-existence 
results of K\"ahler-Einstein metrics to general Einstein metrics. For instance: 

\begin{theorem}\label{t4.3} \cite{L1}. Suppose $(M, g)$ is Einstein, 
non-flat, and $M$ admits an almost complex structure $J$. With respect to 
the orientation and $spin^{c}$ structure induced by $J$, suppose the mod 
2 Seiberg-Witten invariant $\eta_{c}(M) \neq 0$. Then
\begin{equation} \label{e4.6}
\chi (M) \geq  3\tau (M), 
\end{equation}
with equality if and only if $(M, g)$ is homothetic to a complex 
hyperbolic space-form ${\mathbb C}{\mathbb H}^{2}/\Gamma$.
\end{theorem}
\begin{corollary}\label{c4.4} \cite{L1}. The locally symmetric metric $g_{0}$ 
on any compact space-form $M = {\mathbb C}{\mathbb H}^{2}/\Gamma$ is the unique 
Einstein metric on $M$, up to rescaling and isometry.
\end{corollary}

  {\bf Proof:} By a result of Witten and Kronheimer, the manifold $M$ has 
non-zero Seiberg-Witten invariant. Theorem 4.3 then implies that any Einstein 
metric on $M$ is homothetic to a locally symmetric metric, and the result 
then follows from Mostow rigidity.
{\endproof}

  Next we turn to a very different perspective on rigidity and 
non-existence results, developed by Besson-Courtois-Gallot \cite{BCG}. 
This holds in all dimensions, so in the following we assume $M = M^{n}$ 
is a compact $n$-dimensional manifold. 

  For a given compact Riemannian manifold $(M, g)$ let $\widetilde M$ 
denote the universal cover and define the volume entropy of $(M, g)$ by
\begin{equation} \label{e4.7}
h(M, g) = \lim_{r\rightarrow\infty}\frac{\ln vol B_{x}(r)}{r}, 
\end{equation}
where $B_{x}(r)$ is the geodesic $r$-ball about $x \in \widetilde M$. It is 
easy to see this is independent of $x$. Note also that $h(M, g)$ scales 
inversely to the distance, so that $h^{n}(M, g)vol_{g}M$ is scale-invariant. 
If $(M, g_{0})$ is hyperbolic, then a simple computation gives 
$h(M, g_{0}) = n-1$. More generally, if $(M, g)$ is Einstein, then by the 
Bishop-Gromov volume comparison theorem, 
\begin{equation}\label{e4.8}
h(M, g) \leq  (n-1)\sqrt{\frac{\lambda_{-}}{n-1}},
\end{equation}
where $\lambda_{-} = -\min(0, \lambda)$, with equality if and only if 
$(M, g)$ is hyperbolic. 

\begin{theorem}\label{t4.5} \cite{BCG}. Let $(X, g_{0})$ be a compact oriented 
locally symmetric space of negative curvature, and let $M$ be any compact manifold 
with $dim\,M = dim\,X$. Suppose $f: M \rightarrow X$ is smooth. Then for any 
metric $g$ on $M$, one has
\begin{equation}\label{e4.9}
h^{n}(M, g)vol_{g}(M) \geq  |deg f|h^{n}(X, g_{0})vol_{g_{0}}(X),
\end{equation}
with equality if and only if $f$ is a covering map and $g$ is locally 
homothetic to $g_{0}$. 
\end{theorem}

{\bf Sketch of Proof:} The main ideas are already present in the simplest case, 
where $X = M$ is hyperbolic, of curvature -1, and so we assume this in the 
following. Let $L^{2}(S^{n-1}(\infty), d\theta)$ be the Hilbert space of 
$L^{2}$ functions on the sphere at infinity $S^{n-1}(\infty)$ of the hyperbolic 
space ${\mathbb H}^{n}(-1) = \widetilde X$. Let $S_{+}^{\infty}$ be the space 
of positive functions of norm 1 in $L^{2}(S^{n-1}(\infty), d\theta)$. The 
central objects of study are $\pi_{1}(M)$-equivariant Lipschitz maps
\begin{equation}\label{e4.10}
\Phi: \widetilde M \rightarrow  S_{+}^{\infty}.
\end{equation}
The canonical example here is the square root of the Poisson kernel,
$$\Phi_{0} = \sqrt{p_{0}},$$
where $p_{0}(x, \theta) = e^{-(n-1)\beta_{\theta}(x)}$ and $\beta_{\theta}$ is 
the Busemann function associated with the base point $\theta\in S^{n-1}(\infty)$. 
Moreover, for any metric $g$ on $M$, one has another natural class of examples, 
given by $\Phi_{c}(x,\theta) = \Psi_{c}(x,\theta )/|\Psi_{c}(x,\theta)|_{L^{2}(d\theta)}$, 
where
$$\Psi_{c}(x,\theta ) = (\int_{M}e^{-cd(x,y)}p_{0}(y,\theta)dv_{g}(y))^{1/2}.$$
Here $c$ is any number satisfying $c > h(M, g)$, so that $\Psi_{c}(x, \cdot)$ is 
well-defined in $L^{2}(S^{n-1}(\infty), d\theta)$. 
 
   In general, any such $\Phi$ induces a (possibly degenerate) metric 
$g_{\Phi}$ on $M$ by pullback, i.e.
$$g_{\Phi} = \Phi^{*}(g_{can}), $$
where $g_{can}$ is the canonical $L^{2}$ metric (product of $L^{2}$ 
functions) on $S_{+}^{\infty}$. A simple computation shows that
$$g_{\sqrt{p_{0}}} = \frac{(n-1)^{2}}{4n}g_{0},$$
so that the embedding by the Poisson kernel is a homothety. Further 
straightforward computation shows that 
$$g_{\Phi_{c}} \leq \frac{c^{2}}{4}g.$$

 Now define the spherical volume $vol_{Sph}(M)$ to be $\inf vol_{g_{\Phi}}(M)$, 
where the $\inf$ is taken over all $\Phi$ as in \eqref{e4.10}. A simple analysis 
using the functions $\Phi_{c}$ gives an upper bound:
\begin{equation}\label{e4.11}
vol_{Sph}(M) \leq  (\frac{h(g)^{2}}{4n})^{n/2}vol(M, g),
\end{equation}
for any metric $g$ on $M$. In the case at hand, the result \eqref{e4.9} then 
follows from the claim that 
\begin{equation}\label{e4.12}
vol_{Sph}(M) =  (\frac{h(g_{0})^{2}}{4n})^{n/2}vol(M, g_{0}).
\end{equation}
This is proved by exhibiting a calibration form for the Poisson kernel 
embedding. In slightly more detail, for any absolutely continuous measure 
$d\mu$ on $S^{n-1}(\infty)$, define the barycenter $B(\mu)$ to be the 
unique $x\in \widetilde X$ such that
$$\int_{S(\infty )}d\beta_{(x,\theta )}(v)d\mu (\theta ) = 0, $$
for any $v\in T_{x}\widetilde X$. Standard convexity arguments in hyperbolic 
geometry show that $x$ is uniquely defined. This defines a $\Gamma$-equivariant 
map
$$\pi : S_{+}^{\infty} \rightarrow  \widetilde X, \ \ \phi  \rightarrow  
B(\phi^{2}(\theta )d\theta). $$
Now let $\omega_{0}$ be the volume form of the constant curvature metric 
$g_{0}$ on $\widetilde X$. One then shows that the closed $n$-form 
$\pi^{*}\omega_{0}$ on $S_{+}^{\infty}$ is a calibration for the 
Poisson kernel embedding $x \rightarrow  \sqrt{p_{0}(x)}$ of comass 
$(\frac{4n}{h_{0}^{2}})^{n/2}$, which gives \eqref{e4.12}. 

{\endproof}

 Theorem 4.5 gives easily the following uniqueness or rigidity result 
for Einstein metrics on hyperbolic 4-manifolds; note in particular that 
this result gives a new proof of the Mostow rigidity theorem for 
hyperbolic metrics. 

\begin{corollary}\label{c4.6} \cite{BCG}. Suppose $N$ is a compact manifold homotopy 
equivalent to a hyperbolic 4-manifold $(M, g_{0})$. Then $N$ admits an Einstein 
metric only if $N$ is diffeomorphic to $M$, and moreover, $g_{0}$ is the unique 
Einstein metric on $N = M$, up to scaling and isometry.
\end{corollary}

{\bf Proof:} By \eqref{e4.2} and \eqref{e4.3}, for $(N, g)$ Einstein,
$$2\chi \pm 3\tau  = \frac{1}{4\pi^{2}}\int_{N} 2|W_{\pm}|^{2} + \frac{s^{2}}{24} 
- \frac{|z|^{2}}{2} \geq  \frac{1}{6\pi^{2}}\lambda^{2}vol(N, g).$$
For the hyperbolic metric, this gives $2\chi \pm 3\tau  = \frac{3}{2\pi^{2}}
vol(M, g_{0})$, and hence
$$vol(M, g_{0}) \geq   \frac{\lambda^{2}}{9}vol(N, g).$$
Combining this with Theorem 4.5 and \eqref{e4.8} gives
$$3^{4}(\frac{\lambda_{-}}{3})^{2}vol(N, g) \geq  h^{n}(N, g)vol(N, g) 
\geq  3^{4}vol(M, g_{0}) \geq  3^{4}(\frac{\lambda}{3})^{2}vol(N, g). $$
It follows that $\lambda < 0$ and all the inequalities above are 
equalities. The rest of the proof follows from the rigidity statement 
in Theorem 4.5. 
{\endproof}

  In contrast to the Ricci-flat and Ricci-negative results mentioned above, 
there are currently no rigidity or uniqueness results for Ricci-positive 
Einstein metrics, for instance for the standard metrics on $S^{4}$ or 
${\mathbb C}{\mathbb P}^{2}$. It is known, cf.~\cite{Bes}, that the standard 
metrics on these spaces are locally rigid, i.e.~the metrics are isolated 
points in the moduli space of Einstein metrics. 

   An extension of the reasoning in Corollary 4.6 gives the following 
non-existence result.

\begin{corollary}\label{c4.7} \cite{Sa}. For any given values $(k, l)$ with 
$k - l \equiv 0$ $(mod \ 2)$, there are infinitely many non-homeomorphic closed 
4-manifolds $X_{i}$ which satisfy $(\chi (X_{i}), \tau (X_{i})) = (k, l)$, 
and which admit no Einstein metric. 
\end{corollary}

{\bf Sketch of Proof:} The idea is to take connected sums of a hyperbolic manifold 
$M$ as above with copies of $\pm{\mathbb C}{\mathbb P}^{2}$, $S^{2}\times S^{2}$ or 
$S^{2}\times T^{2}$ and use the degree theory part of Theorem 4.5. 
{\endproof}

  There are a several further interesting obstructions to the existence of 
Einstein metrics on 4-manifolds, but for lack of space we will forgo a detailed 
discussion. First, there is an improvement of the Hitchin-Thorpe inequality 
due to Gromov \cite{Gr2}, when the manifold $M$ has non-zero simplicial volume, 
cf.~in particular \cite{Kt1}, \cite{Kt4}. Next, based on information from the Seiberg-Witten 
invariants developed by LeBrun in \cite{L2}, Kotschick, LeBrun and many others have 
found a wide variety of simply connected 4-manifolds of a fixed homeomorphism 
type, which have an Einstein metric for one smooth structure, but without Einstein 
metrics for other smooth structures, analogous to the discussion following 
Theorem 4.2. Also, in the same context, there are Einstein metrics with 
$\lambda > 0$ for one smooth structure and with $\lambda < 0$ for a different 
smooth structure; we refer to \cite{L2,L3}, \cite{Kt2,Kt3} and references therein 
for further details.

\section{Moduli spaces I.}
\setcounter{equation}{0}

 In this section, we discuss various aspects of the moduli space ${\mathcal E}$ 
of Einstein metrics on a given compact 4-manifold $M$. We begin with local results, 
(which hold in all dimensions), and then pass to more global issues on the 
structure of ${\mathcal E}$. 

 The Einstein equations \eqref{e1.1} are invariant under scaling, and so 
throughout \S 5 and \S 6, we assume all metrics are normalized to have 
unit volume. Let ${\mathbb E} = {\mathbb E}(M)$ denote the space of all 
(unit volume) Einstein metrics on a given manifold $M$, viewed as a subset of 
the space $Met(M)$ of all unit volume Riemannian metrics on $M$. As noted in 
\S 2, Einstein metrics are $C^{\infty}$ smooth, in fact real-analytic, in suitable 
local coordinate systems. The group ${\mathcal D}$ of $C^{\infty}$ diffeomorphisms 
acts continuously on ${\mathbb E}$ and the quotient
\begin{equation}\label{e5.1}
{\mathcal E} = {\mathcal E}(M)
\end{equation}
is the moduli space of Einstein metrics on $M$. It is standard that 
${\mathcal E}$ is Hausdorff, with countably many components, cf.~\cite{Bes}. 
As noted in \S 2, Einstein metrics are critical points of the total scalar 
curvature \eqref{e2.2}. Hence $s_{g}$, or equivalently $\lambda$, is 
constant on each component of ${\mathcal E}$.

  The Einstein equations \eqref{e1.1} are not elliptic, due to their invariance 
under diffeomorphisms. In fact if $\hat E(g) = Ric_{g} - \lambda g$ denotes the 
Einstein operator, then the linearization $D\hat E$ is given by
$$2D\hat E_{g}(h) = D^{*}Dh - 2R(h) - 2\delta^{*}\beta (h), $$
where $R(h)$ is the action of the curvature tensor on the space $S^{2}(M)$ 
of symmetric bilinear forms on $M$, $\beta$ is the Bianchi operator, 
$\beta (h) = \delta h + \frac{1}{2}dtr h$ and $D^{*}D$ is the rough Laplacian. 
One has $D\hat E_{g}(\delta^{*}X) = 0$, for any vector field $X$ on $M$, so 
that $Ker\,D\hat E_{g}$ is infinite dimensional. 

  As is usual in geometrically covariant problems, one needs to fix a 
gauge transverse to the action of ${\mathcal D}$ to obtain an elliptic 
system. To do this, first pass to the usual Einstein operator used in physics, 
\begin{equation}\label{e5.2}
E(g) = Ric_{g} - \frac{s}{2}g + \Lambda g,
\end{equation}
where $\Lambda = \frac{n-2}{2}\lambda$, $n = dim M$. Then $E^{-1}(0)$ consists 
of Eintein metrics \eqref{e1.1}. Choose a background metric $g_{0}\in {\mathbb E}$ 
and consider the divergence-gauged Einstein operator
\begin{equation}\label{e5.3}
\Phi_{g_{0}}(g) = Ric_{g} - \frac{s}{2}g + \Lambda g + 
\delta_{g}^{*}\delta_{g_{0}}(g),
\end{equation}
where $\delta_{g_{0}}$ is the divergence operator with respect to $g_{0}$. 
The linearization of $\Phi$ at $g = g_{0}$ is given by 
\begin{equation}\label{e5.4}
2(D\Phi)_{g_{0}}(h) = L(h) = D^{*}Dh - 2R(h) - D^{2}tr h - \delta \delta h\,g 
+ \Delta tr h \,g + \frac{s}{n}tr h\,g .
\end{equation}
The operator $L$ is self-adjoint, and elliptic when $n \geq 3$. Thus $L$ is Fredholm, 
and so has finite dimensional kernel and cokernel with closed range. Let ${\mathbb Z} 
= \Phi^{-1}(0)$ be the zero-set of $\Phi$. We observe that for metrics $g$
close to $g_{0}$ one has 
\begin{equation}\label{e5.5}
{\mathbb Z}  \subset  {\mathbb E}.
\end{equation}
To see this, apply $\delta_{g}$ to \eqref{e5.3}. By the Bianchi identity, 
$\delta_{g}E(g) = 0$, and hence, for $g \in {\mathbb Z}$ and $V = \delta_{g_{0}}g$, 
one has
$$\delta \delta^{*}V = 0.$$
Pairing this with $V$ and integrating by parts gives $\delta^{*}V = 0$ which implies 
\eqref{e5.5}. 

  On the other hand, via the Ebin slice theorem, for any $g\in {\mathbb E}$ 
near $g_{0}$ there exists $\phi\in{\mathcal D}$ such that $\phi^{*}g$ is in 
divergence-free gauge, i.e.
$$\delta_{g_{0}}(\phi^{*}g) = 0. $$
Thus ${\mathbb Z}$ is a local slice for ${\mathbb E}$, transverse to the action 
of ${\mathcal D}$. This leads to the following result of Koiso. 

\begin{theorem} \label{t5.1} \cite{Koi}. Near any $g_{0}\in {\mathbb E}$, the space 
${\mathbb Z} \subset {\mathbb E}$ has the structure of a finite dimensional 
real-analytic subvariety of a smooth, finite dimensional manifold ${\mathbb W}$. 
Further the tangent space $T_{g_{0}}{\mathbb W}$ consists of the space of 
essential infinitesimal Einstein deformations of $g_{0}$. 
\end{theorem}

{\bf Proof:} Let $\pi: S^{2}(M) \rightarrow Im\, L$ be the 
projection onto Im $L = Im \,D\Phi$, and recall that $Im\,L$ is of 
finite codimension in $S^{2}(M)$. The composition $F = \pi\circ\Phi: 
Met(M) \rightarrow Im\,L$ is then a submersion near $g_{0}$, 
i.e.~the derivative $DF$ is surjective, with splitting kernel. It follows 
from the implicit function theorem in Banach spaces that ${\mathbb W} \equiv 
F^{-1}(0)$ is smooth submanifold of $Met(M)$ near $g_{0}$ of dimension 
$dim \,Ker\, L = dim\, Coker\, L < \infty$. 

  The tangent space to ${\mathbb W}$ at $g = g_{0}$ equals $Ker\, L$. 
For $h \in Ker\,L$, the same arguments establishing \eqref{e5.5} show that 
$\delta h = 0$. Using this, and taking the trace of \eqref{e5.4} then gives
$$\Delta tr h + {\tfrac{s}{n}}tr h = 0.$$
If $s \leq 0$, then it is immediate that $tr h = 0$. If $s > 0$, this 
conclusion also follows from the Lichnerowicz estimate on the first 
eigenvalue of the Laplacian: $\lambda_{1} \geq \frac{s}{n-1}$. Hence, 
$Ker\,L$ consists of the forms $h$ satisfying
\begin{equation}\label{e5.6}
D^{*}Dh - 2R(h) = 0.
\end{equation}
This, together with the conditions $\delta h = tr h = 0$ are the 
equations for essential infinitesimal Einstein deformations. 

  Finally, $\Phi$ is a real-analytic function on ${\mathbb W}$, so the 
zero set ${\mathbb Z}$ is a real-analytic subvariety of ${\mathbb W}$. 
{\endproof}

\begin{remark}\label{r5.2} {\bf (i).} 
{\rm In \cite{Koi}, Koiso finds examples where one has a strict inclusion 
$${\mathbb Z} \subset {\mathbb W} , $$
and hence there are situations where infinitesimal Einstein 
deformations are not tangent to a curve of Einstein metrics in 
${\mathbb E}$; the moduli space is non-integrable or obstructed. 
This occurs for instance on ${\mathbb C}{\mathbb P}^{1}\times 
{\mathbb C}{\mathbb P}^{2k}$. However, there are no examples 
where this occurs in dimension 4. 

 In fact, there are still no examples where ${\mathbb Z}$ is not a 
finite dimensional manifold. Is the non-integrability related to the existence 
of Killing fields? For instance, if there are no Killing fields on $(M, g)$, 
(for example $\lambda < 0$), is ${\mathbb Z}  = {\mathbb W}$ near $g$? 

 {\bf (ii).} Is there any method to compute the dimension of components 
of ${\mathcal E}$? 

{\bf (iii).} There are several well-known local rigidity results for 
Einstein metrics under various curvature conditions, which show that 
a given metric $g \in {\mathcal E}$ is an isolated point in the full 
moduli space. This is the case for instance for Einstein metrics of 
strictly negative sectional curvature, or for irreducible symmetric spaces, 
cf.~\cite{Bes} for further discussion. All of these results follow from an 
analysis (generally algebraic) on solutions of \eqref{e5.6}. }
\end{remark}

 The first global results on the structure of moduli space ${\mathcal E}$ 
of unit volume Einstein metrics on a given 4-manifold $M$ were obtained in 
\cite{An1}, \cite{BKN}, \cite{N1}. These results bear some similarities to 
the results of Uhlenbeck \cite{FU} on the moduli space of self-dual Yang-Mills 
fields. 

 To describe the situation, we first need the following definition. An Einstein 
orbifold $(V, g)$ associated to a 4-manifold $M$ is a 4-dimensional orbifold, with 
a finite number of singular points $q_{k}$, each having a neighborhood homeomorphic 
to the cone $C(S^{3}/\Gamma)$, where $\Gamma \neq \{e\}$ is a finite subgroup 
of $SO(4)$. Let $V_{0} = V\setminus \cup q_{k}$ be the regular, (smooth manifold), 
set of $V$. Then $g$ is a smooth Einstein metric on $V_{0}$, which extends smoothly 
over $\{q_{k}\}$ in local finite covers. The manifold $M$ is a resolution of 
$V$ in the sense that there is a continuous surjection $\pi : M 
\rightarrow  V$ such that $\pi|_{\pi^{-1}(V_{0})}: \pi^{-1}(V_{0}) 
\rightarrow  V_{0}$ is a diffeomorphism onto $V_{0}$. In particular, 
$V$ is compact.  

\begin{theorem}\label{t5.3} \cite{An1}, \cite{BKN}, \cite{N1}. The 
completion $\bar{\mathcal E}_{GH}$ of ${\mathcal E}$ in the Gromov-Hausdorff 
topology consists of ${\mathcal E}$ together with unit volume Einstein orbifold 
metrics associated to $M$. 

 Moreover, the completion is locally compact, in that any sequence 
$g_{i} \in \bar{\mathcal E}_{GH}$, bounded in the Gromov-Hausdorff topology, 
has a subsequence converging to an Einstein orbifold associated to $M$. 
There is a uniform bound on the number of orbifold singularities and the 
order of the local groups $\Gamma$ in terms of $\chi (M)$.
\end{theorem}

{\bf Sketch of Proof:} By \eqref{e4.2}, unit volume Einstein metrics on $M$ have 
a lower bound on their scalar curvature and hence a uniform lower bound on their 
Ricci curvature. The completion in the Gromov-Hausdorff topology is then equivalent 
to the completion with respect to a diameter bound, so that any Cauchy sequence 
$\{g_{i}\} \in {\mathbb E}$ satisfies
\begin{equation} \label{e5.7}
vol_{g_{i}}M = 1, \  s_{g_{i}} \geq s_{0} > 0, \ diam_{g_{i}}M \leq  D,  
\end{equation}
for some $D = D(g_{i}) < \infty$. Gromov's weak compactness theorem \cite{Gr1} 
implies that the Cauchy sequence $\{g_{i}\}$ converges in the Gromov-Hausdorff 
topology to a complete length space $(X, d_{\infty})$. One needs then to 
understand the structure of the limit $(X, d_{\infty})$. 

 Given $r > 0$, let $\{x_{k}\}$ be a maximal $r/2$ separated set, (depending on $i$), 
in $(M, g_{i})$. Thus, the geodesic balls $B_{x_{k}}(\frac{r}{2})$ are disjoint 
while the balls $B_{x_{k}}(r)$ cover $M$. Choose a fixed $\delta_{0} > 0$ small, 
and let
\begin{equation} \label{e5.8}
G_{i}^{r} = \cup\{B_{x_{k}}(r): \int_{B_{x_{k}}(2r)}|R|^{2}dV < 
\delta_{0}\}, 
\end{equation}
and similarly, let
\begin{equation} \label{e5.9}
B_{i}^{r} = \cup\{B_{x_{k}}(r): \int_{B_{x_{k}}(2r)}|R|^{2}dV \geq  
\delta_{0}\}. 
\end{equation}
All quantities here are with respect to $(M_{i}, g_{i})$. For each $i$, one has 
$M_{i} = G_{i}^{r}\cup B_{i}^{r}$. Observe via \eqref{e4.2} that there is uniform 
bound $K = 8\pi^{2}\chi(M)$ on the number $Q_{i}^{r}$ of $r$-balls in $B_{i}^{r}$:
\begin{equation} \label{e5.10}
Q_{i}^{r} \leq  \frac{K}{\delta_{0}}.
\end{equation}

  Einstein metrics satisfy the inequality
$$\Delta |R| + c|R|^{2} \geq 0,$$
where $c$ is a constant, depending only on dimension. Using this together with the 
deGiorgi-Nash-Moser method in elliptic PDE, cf.~\cite{GT}, one shows that for 
$\delta_{0}$ small, depending only on $D$ in \eqref{e5.7}, one 
has the $L^{\infty}$ estimate
\begin{equation} \label{e5.11}
|R| \leq  C\delta_{0}r^{-2} \ \ {\rm on} \ \  G_{i}^{r}. 
\end{equation}
In fact \eqref{e5.11} holds on the $r/2$ thickening of $G_{i}^{r}$. It then 
follows from the smooth Gromov compactness theorem that for any given 
$r > 0$, a subsequence of $G_{i}^{r}$ converges in the $C^{1,\alpha}$ 
topology to a limit manifold $G_{\infty}^{r}$ with limit $C^{1,\alpha}$ 
metric $g_{\infty}^{r}$. In particular, $G_{\infty}^{r}$ and $G_{i}^{r}$ are 
diffeomorphic, for $i$ large, and there exist smooth embeddings 
$F_{i}^{r}: G_{\infty}^{r} \rightarrow  G_{i}^{r} \subset  M_{i}$ such 
that $(F_{i}^{r})^{*}(g_{i})$ converges in $C^{1,\alpha}$ to 
$g_{\infty}^{r}$. Via regularity of the Einstein equation as in \eqref{e2.2}, 
both the limit metric $g_{\infty}^{r}$ and the convergence are in fact 
$C^{\infty}$ smooth. 

 Now choose a sequence $r_{j} \rightarrow 0$, with $r_{j+1} = 
\frac{1}{2}r_{j},$ and perform the above construction for each $j$. Let 
$G_{i}(r_{m}) = \{x\in (M_{i}, g_{i}): x\in G_{i}^{j}, {\rm for \ some} 
\ j \leq  m\}$, so that one has inclusions
$$G_{i}(r_{1}) \subset G_{i}(r_{2}) \subset  ... \subset  M_{i} $$
By the argument above, each $G_{i}(r_{m}) \subset  (M_{i}, g_{i}),$ for 
$m$ fixed, has a subsequence converging smoothly to a limit 
$G_{\infty}(r_{m})$. Clearly $G_{\infty}(r_{m}) \subset G_{\infty}(r_{m+1})$ 
and we set
$$G = \cup_{1}^{\infty}G(r_{m}),$$
with the induced metric $g_{\infty}$. Thus, $(G_{\infty}, g_{\infty})$ 
is $C^{\infty}$ smooth and for any $m$, there are smooth 
embeddings $F_{i}^{m}: G(r_{m}) \rightarrow  M_{i}$, for $i$ 
sufficiently large, such that $(F_{i}^{m})^{*}(g_{i})$ converges smoothly 
to the metric $g_{\infty}$. 

 Let $\bar G$ be the metric completion of $G$ with respect to $g_{\infty}$. 
Then there is a finite set of points $q_{k}$, $k = 1, \cdots , Q$, such that
$$\bar G = G \cup  \{q_{k}\}.$$
This follows since there is a uniform upper bound \eqref{e5.10} on the 
cardinality of $B_{i}^{r},$ for all $r$ small, and all $i$, independent 
of $r, i$. It is then easy to see that a subsequence of $(M_{i}, g_{i})$ 
converges to the length space $(\bar G, g_{\infty})$ in the 
Gromov-Hausdorff topology, so that $X = \bar G$. 

 It remains to prove that $\bar G$ is an orbifold, with orbifold 
singular points $\{q_{k}\}$, i.e.~$\bar G = V$. This follows by an 
analysis of the tangent cone of the limit metric $g_{\infty}$ near each 
$q_{k}$, i.e.~by a blow-up analysis. The curvature of $G$ is locally 
bounded in $L^{p}$, for any $p < \infty$. Further, by lower 
semi-continuity of the norm under weak convergence, the $L^{2}$ norm 
of the curvature on $(G, g_{\infty})$ is globally bounded:
\begin{equation} \label{e5.12}
\int_{G}|R|^{2}dV_{g_{\infty}} \leq  8\pi^{2}\chi(M) . 
\end{equation}
In particular, for any $q = q_{k}\in \bar G,$
\begin{equation} \label{e5.13}
\int_{A_{q}(\frac{1}{2}r,2r)}|R|^{2}dV \rightarrow  0, \ \ {\rm as} \ \  
r \rightarrow  0, 
\end{equation}
and hence, as in \eqref{e5.11}, near any singular point $q$ one has
\begin{equation} \label{e5.14}
|R| \leq \varepsilon(r)r^{-2},
\end{equation}
where $r(x) = dist(x, q)$ and $\varepsilon(r) \rightarrow 0$ as 
$r \rightarrow 0$. 

 Let $s_{j} = 2^{-j},$ for $j$ large, and rescale the metric 
$g_{\infty}$ on $G$ near $q$ by $s_{j}^{-2}$, i.e.~consider the metrics 
$\bar g_{j} = s_{j}^{-2}\cdot  g_{\infty}$. The bounds \eqref{e5.7}, 
together with the Bishop-Gromov volume comparison theorem and the 
curvature bound \eqref{e5.14} imply that a subsequence of 
$(G\setminus q, \bar g_{j}, q)$ converges, modulo diffeomorphisms, 
smoothly to a flat limit $(T^{\infty}, \bar g^{\infty})$. 
Next one shows that $T^{\infty}$ has a bounded number of components and 
the metric completion $\bar T^{\infty}$ of each component of $T^{\infty}$ 
has a single isolated singularity \{0\}. Thus, $\bar T^{\infty}$ is a 
finite collection of complete flat manifolds joined at a single 
isolated singularity $\{0\}$. From this, it follows easily that 
$\bar T^{\infty}$ is isometric to a union of flat cones 
$C(S^{n-1}/\Gamma_{j})$. By the smooth convergence, this (unique) 
structure on the limit is equivalent to the structure of $(G, g_{\infty})$ 
on small scales near the singular point $q$. Via the Cheeger-Gromoll 
splitting theorem, one proves that the regular set $V_{0}$ is locally 
connected and hence each singular point is an orbifold singularity. 
It is proved in \cite{BKN} that in the local finite cover resolving a 
singularity $q$, the lifted metric $g_{\infty}$ extends smoothly across 
the origin. 

 The remaining parts of Theorem 1.1 are now easily established, via 
\eqref{e5.7}, \eqref{e4.2} and the volume comparison theorem. For the 
proof of Theorem 5.3 in the case that $g_{i}$ are Einstein orbifold 
metrics associated to $M$, see \cite{An3}. 
{\endproof}

  A main point of the Uhlenbeck completion for 
self-dual connections is that the completion is {\it compact}. In the current 
context of Einstein metrics, consider the components ${\mathcal E}_{\lambda_{0}}$ 
of ${\mathcal M}$ for which 
\begin{equation} \label{e5.15}
\lambda \geq  (n-1)\lambda_{0} >  0. 
\end{equation}
Myers' theorem, (for manifolds of positive Ricci curvature), then implies that
$$diam_{g}M \leq \pi / \sqrt{\lambda_{0}},$$
so that \eqref{e5.7} holds automatically. Hence, the completion of 
${\mathcal E}_{\lambda_{0}}$ in the Gromov-Hausdorff topology is compact. 
However, this is certainly not the case when $\lambda \leq 0$. For example, 
the moduli space ${\mathcal E}$ or $\bar{\mathcal E}_{GH}$ on a torus $T^{4}$, 
or on a product $\Sigma_{g_{1}}\times \Sigma_{g_{2}}$ of surfaces of genus 
at least 2, is certainly not compact. Thus, for a better understanding one 
needs to consider what happens when the Gromov-Hausdorff distance goes to 
infinity. This will be discussed further in \S 6. 

\medskip

 There is another strong difference compared with the Uhlenbeck 
completion. Namely, the frontier $\partial_{o}{\mathcal E} = 
\bar{\mathcal E}_{GH}\setminus {\mathcal E}$ should not really be thought 
of as a boundary, but instead as a filling in of ``missing pieces'' in 
${\mathcal E}$. For example, as discussed in \S 6, in the case of K3 surfaces 
this frontier consists of subvarieties of codimension 3 in $\bar{\mathcal E}$ 
and so does not form a boundary in the sense of a wall at which the moduli 
space comes to an end. 

  Although there is currently very little evidence, we venture the following 
(optimistic) conjecture, which would confirm this picture in general:

\noindent
{\bf Conjecture.} The space $\partial_{o}{\mathcal E} \subset \bar{\mathcal E}_{GH}$ 
of Einstein orbifold metrics associated to $M$ is of codimension at least 2 in 
$\bar{\mathcal E}_{GH}$.

\medskip

 The orbifold limits $(V, g_{\infty})$ arise from the ``bubbling off'' of 
so-called gravitational instantons, (again in analogy to the case of Yang-Mills 
fields). These spaces, called EALE spaces here, are complete Ricci-flat metrics 
$(N, g)$ which have curvature in $L^{2}$, 
\begin{equation}\label{e5.16}
\int_{N}|R|^{2}dV_{g} < \infty,
\end{equation}
and which are ALE, (asymptoticaly locally Euclidean), in that the metric 
$g$ at infinity is asymptotic to a flat cone $C(S^{3}/\Gamma)$, 
where $\{e\} \neq \Gamma \subset SO(4)$ is a finite subgroup. Thus, 
outside a compact set $K \subset N$, there is a finite cover of 
$N\setminus K$ which is diffeomorphic to ${\mathbb R}^{4}\setminus B$, 
and a chart in which the (lifted) metric $g$ has the form
$$|g_{ij} - \delta_{ij}| = \varepsilon (r), \ \ |R| = \varepsilon (r)r^{-2}, $$
where $\varepsilon(r) \rightarrow 0$ as $r \rightarrow \infty$. It is proved 
in \cite{BKN} that in fact one has $\varepsilon(r) = r^{-2}$. A similar definition 
holds for EALE orbifolds. 

 To describe this bubbling process, let $q$ be a singular point of the limit 
$V$. If $x_{i}\in (M, g_{i})$ is any sequence of points such that $x_{i} 
\rightarrow  q$ in the Gromov-Hausdorff topology, then the curvature of 
$(M, g_{i})$ blows up near $x_{i}$, i.e.~diverges to infinity, as 
$i \rightarrow \infty$. If one rescales the metrics $g_{i}$ so that the 
curvature remains bounded near $x_{i}$, then a subsequence converges to 
a complete EALE space $(N, g)$. Blowing this limit $(N, g)$ down, 
i.e.~rescaling $g$ by factors converging to 0, gives a spherical cone 
at a single vertex $\{0\}$.

 However, the curvature of $(M, g_{i})$ may diverge to infinity at a number 
of different scales near any singular point $q$, giving rise to a 
collection of such EALE spaces associated with each scale. This gives 
rise to a so-called "bubble-tree" of EALE spaces and scales. The 
structure of the limit orbifold $(V, g_{\infty})$ near any singular 
point $q$ is recaptured by the structure at infinity of the complete 
EALE orbifold corresponding to the {\it  smallest} rate at which the 
curvature of $(M_{i}, g_{i})$ diverges to infinity near $q$; this 
corresponds to the {\it largest} distance scale, since the curvature scale 
corresponds to the inverse square of the distance scale. 

 In more detail, for $x_{i} \rightarrow q$ as above, there is sequence 
of scales $r_{i} = r_{i}^{1} \rightarrow 0$ such that the rescalings 
$(M, r_{i}^{-2}g_{i}, x_{i})$ with $x_{i} \rightarrow q$, converge, 
(in a subsequence), in the pointed Gromov-Hausdorff topology to a 
complete, EALE orbifold $(V^{1}, g^{1})$ with a finite number of 
singular points, and
\begin{equation} \label{e5.17}
\int_{V^{1}}|R|^{2} \geq  \delta_{0}, 
\end{equation}
for a fixed $\delta_{0} > 0$, (as in \eqref{e5.9}). If $(V^{1}, g^{1})$ is not a 
smooth manifold, then there are second level scales $\{r_{i}^{2}\}$ associated 
with each singular point of $V^{1}$; (the scales $r_{i}^{2}$ depend on the 
choice of singular point in $V^{1}$). Rescaling $r_{i}^{-2}g_{i}$ by such 
factors at base points converging to the singular points of $V^{1}$ gives a 
collection of $2^{\rm nd}$ level EALE orbifolds $\{(V^{2}, g^{2})\}$ associated 
with each singular point of $(V^{1}, g^{1})$. Each iteration of this process 
satisfies \eqref{e5.17}, and hence, via \eqref{e4.2}, (as in \eqref{e5.10}), 
terminates after a finite number of iterations. At the last stage, corresponding 
to the smallest scales, the resulting blow-up limits are non-flat smooth 
manifolds, cf.~\cite{Ba}, \cite{AC} for full details. 

   The topology of the original manifold $M$ may then be reconstituted from 
that of $V_{0}$, and the scale of orbifolds associated with each singular 
point $q\in V$. In particular, the homology groups of $M$ are determined, 
(by the Mayer-Vietoris sequence), by the homology of $V$ and the homology 
of the collection of EALE orbifold spaces $\{V^{j}\}$. Thus the orbifold 
singularities correspond to a generalized connected sum decomposition of 
$M$, in that $M$ is the union of the regular set $V_{0}$ with a finite 
collection of EALE spaces:
\begin{equation}\label{e5.18}
M = V_{0}\cup\{N_{m}\},
\end{equation}
where the union is along non-trivial spherical space forms. The collection 
$\{N_{m}\}$ itself could well consist of orbifolds, in which case it also 
splits inductively as a union along spherical space forms.  

\medskip

 Clearly then it is important to understand the geometry and topology of 
these EALE spaces in detail. In the special case where $(N, g)$ is 
simply connected and K\"ahler, (and so hyperk\"ahler), one has a complete 
description and classification due to Kronheimer. To describe this, let 
$\Gamma \subset SU(2)$ be a finite group. Then $\Gamma$ belongs to 
one of the following five classes.

 $\bullet$ $A_{k}$: $\Gamma$ is the cyclic group ${\mathbb Z}_{k+1}.$

 $\bullet$ $D_{k}$: $\Gamma$ is the binary dihedral group of order $4(k-2).$

 $\bullet$ $E_{6}$, $E_{7}$, $E_{8}$: $\Gamma$ is the binary tetrahedral, 
binary octahedral group or binary icosahedral group respectively. 

 The quotient ${\mathbb C}^{2}/\Gamma$ is called a rational double point. 
Let $\pi : N_{\Gamma} \rightarrow{\mathbb C}^{2}/\Gamma$ be a minimal 
resolution of ${\mathbb C}^{2}/\Gamma$. The exceptional divisor $E = 
\pi^{-1}(0)$ is a union of ${\mathbb C}{\mathbb P}^{1}$'s, $E = \Sigma_{1} 
+ \cdots + \Sigma_{l}$, and so in particular $H_{2}(N_{\Gamma}, {\mathbb Z}) 
= \oplus_{1}^{l}{\mathbb Z}$ is generated by $\{\Sigma_{i}\}$. Moreover, the 
intersection form $(\Sigma_{i}\cdot \Sigma_{j})$ is given by the negative 
of the Cartan matrix associated to the root system of the Lie group associated 
to $\Gamma$. 

\begin{theorem}\label{t5.4} \cite{Kr1,Kr2}. For any such $N_{\Gamma}$, there is a 
hyperk\"ahler EALE metric $(N_{\Gamma}, g)$, which is uniquely determined by three 
cohomology classes $\alpha_{1}$, $\alpha_{2}$, $\alpha_{3} \in H^{2}(N_{\Gamma}, 
{\mathbb R})$ which satisfy the following nondegeneracy condition: for each 
$\Sigma \in H_{2}(N_{\Gamma}, {\mathbb Z})$ with $\Sigma \cdot \Sigma = -2$, 
there exists $i$ such that $\alpha_{i}(\Sigma) \neq 0$. 

  Conversely, any ALE hyperk\"ahler 4-manifold is diffeomorphic to some 
$N_{\Gamma}$, and any such metric is uniquely determined by the cohomology 
classes of the K\"ahler forms $\alpha_{i}$. 
\end{theorem}

 Kronheimer's result above arose from earlier work of Eguchi-Hanson \cite{EH}, 
Gibbons-Hawking \cite{GH1} and Hitchin \cite{Hi2}, (among others), and it is useful 
to describe in detail the Gibbons-Hawking metrics. These are the 
hyperk\"ahler EALE spaces for which $\Gamma  = A_{k}$. 

 Thus, choose $(k+1)$ points $\{p_{i}\}$ in ${\mathbb R}^{3}$ and let $U = 
{\mathbb R}^{3} \setminus \{p_{i}\}$. Let $\pi_{0}: N_{0} \rightarrow U$ 
be the principal $S^{1}$ bundle whose first Chern class is -1 when restricted 
to a small sphere about any $p_{i}$. Since $H_{2}(U, {\mathbb Z}) = 
\oplus_{1}^{k}{\mathbb Z}$, this uniquely determines the principal $S^{1}$ 
bundle $\pi_{0}$. Moreover, for $r$ small, the domain $\pi_{0}^{-1}(B_{p_{j}}(r)) 
\subset N_{0}$ is diffeomorphic to a punctured 4-ball and so adding a point to 
each such neighborhood gives a closed, non-compact manifold $N$ and a smooth map 
$\pi: N \rightarrow  {\mathbb R}^{3}$ extending the map $\pi_{0}$. If 
$l_{i}$ is a line segment joining $p_{i}$ to $p_{i+1}$, then 
$\pi^{-1}(l_{i})$ is a 2-sphere $S^{2} \subset N$, with self-intersection 
-2. The collection of such 2-spheres generates the homology group 
$H_{2}(N, {\mathbb Z}) = \oplus_{1}^{k}{\mathbb Z}$. Alternately, the 
topology of $N$ is given by plumbing a collection of $k$ spheres $S^{2}$ 
according to the Cartan matrix of $A_{k}$. 

 To describe the Gibbons-Hawking metric, let $V: U \rightarrow {\mathbb R}^{+}$ 
be given by
\begin{equation}\label{e5.19}
V = \frac{1}{2}\sum_{1}^{k}\frac{1}{|p - p_{j}|},
\end{equation}
and note that $V$ is harmonic on $(U, g_{0})$, where $g_{0}$ is the 
Euclidean metric on ${\mathbb R}^{3}$. It is easily verified that the cohomology 
class of the closed 2-form $\frac{1}{2\pi}*dV$ represents the Chern class 
$c_{1}(\pi_{0})$ of the $S^{1}$ bundle $\pi_{0}$ in deRham cohomology. 
Thus, there is a connection form $\omega \in \Omega^{1}(M_{0})$ whose 
curvature is given by $*dV$, so that
$$\pi_{0}^{*}(*dV) = d\omega . $$
The Gibbons-Hawking metric has the simple explicit form
\begin{equation}\label{e5.20}
g_{GH} = V^{-1}w\cdot \omega  + V\pi_{0}^{*}(g_{0}).
\end{equation}
It is easy to see that $g_{GH}$ is smooth, and self-dual, and 
hence (K\"ahler) and Ricci-flat.

 An important special case, corresponding to $k = 2$, is the 
Eguchi-Hanson metric on $TS^{2}$:
\begin{equation}\label{e5.21}
g_{EH} = V^{-1}dr^{2} + V\sigma_{1}^{2} + r^{2}(\sigma_{2}^{2} + 
\sigma_{3}^{2}),
\end{equation}
where $V(r) = (1 - (\frac{a}{r})^{4})$ and $\{\sigma_{i}\}$ is the standard 
coframing of $SO(3) \simeq {\mathbb R}{\mathbb P}^{3}$ . Note that this metric 
has a free, isometric ${\mathbb Z}_{2}$ action, and hence the metric desends to 
a metric on the ${\mathbb Z}_{2}$ quotient, (where $b_{2} = 0$). 

\medskip

  The following is a well-known open problem: are there any non-K\"ahler 
EALE spaces $(N, g)$, (at least if $N$ is simply connected)? In \cite{N2}, 
Nakajima proved that if $N$ is spin, and the group $\Gamma$ is a 
subgroup of $SU(2) \subset SO(4)$, and acting suitably on ${\mathbb C}^{2}$, 
then $N$ is necessarily hyperk\"ahler. On the other hand, since there are 
many Yang-Mills connections which are not self-dual, one might expect 
that the answer to this problem is negative. 

  For later purposes, we note that if one changes the potential function 
$V$ in \eqref{e5.18} to 
\begin{equation}\label{e5.22}
V = 1+ \frac{1}{2}\sum_{1}^{k}\frac{1}{|p - p_{j}|},
\end{equation}
then the metric \eqref{e5.19} remains Ricci-flat and K\"ahler, and gives the 
so-called multi-Taub-NUT metrics, cf.~\cite{GH2}. These are no longer ALE, but 
are ALF, (asymptotically locally flat), with asymptotics of the form 
$S^{3}/\Gamma$, where the $S^{1}$ fiber of the Hopf fibration is shrunk to 
bounded length near infinity. When $k = 1$, the potential \eqref{e5.22} 
gives the (Riemannian) Taub-Nut metric. 

\medskip

 Returning now to the general topological issues as in \eqref{e5.18}, it is clearly 
useful to understand the topology of general non-flat EALE spaces $N$. First, by the 
Cheeger-Gromoll splitting theorem $\partial N$ is connected and so the inclusion 
map $i: \partial N = S^{3}/\Gamma \rightarrow N$ induces a surjection
\begin{equation}\label{e5.23}
\pi_{1}(S^{3}/\Gamma) \rightarrow \pi_{1}(N) \rightarrow 0.
\end{equation}
Next, (by \cite{An1} for instance), one has
\begin{equation}\label{e5.24}
|\pi_{1}(N)| <  |\Gamma|.
\end{equation}
In particular, the universal cover of an EALE space is still EALE. The following 
result essentially appears in \cite{An3}. 

\begin{proposition}\label{p5.5} If $N$ is simply connected, then $N$ has 
the homotopy type of a bouquet of $2$-spheres. In particular,
\begin{equation}\label{e5.25}
b_{2}(N) > 0.
\end{equation}
\end{proposition}

{\bf Proof:} A result of Wu \cite{Wu} applies in this setting, and implies that 
there exists a smooth Morse-Smale exhaustion function $\rho$ on $N$, which 
has non-degenerate critical points of index at most 2; $\rho$ is obtained 
by a smoothing of the distance function from any given point in $N$. All 
flow lines of $\nabla\rho$ either connect critical points of $\rho$, or 
diverge to infinity in $N$, and hence $N$ can be retracted onto the space 
of flow lines connecting critical points. In particular, general $N$, (not 
necessarily simply connected), have the homotopy type of a 2-dimensional 
CW complex. 

  The Morse-Smale condition means that if $x$, $y$ are any distinct critical 
points of $\rho$, with $\rho (y) > \rho(x)$, then the space $M_{x,y}$ of flow 
lines starting at $x$ and ending at $y$ is a smooth submanifold of $N$, 
(given by the intersection of the unstable manifold of $x$ and the stable 
manifold of $y$). Further, 
$$dim M_{x,y} = ind_{y} - ind_{x}, $$
where $ind_{x}$ is the index of the critical point $x$. 

 Thus, for any $x \neq y$, one has dim $M_{x,y} = 1$ or 2. There are a finite 
number of 1-dimensional components $M_{x,y}$, given by arcs from $x$ to $y$. 
The closure $\bar M_{x,y}$ of a 2-dimensional component is either a 2-sphere 
$S^{2}$, or a disc $D^{2}$ with edge identifications along $S^{1} = \partial 
D^{2}$, cf.~\cite{We} for instance. The components of the latter type have 
non-trivial $\pi_{1}$, and the result then follows from the Seifert-van Kampen 
theorem for the fundamental group. 
{\endproof}

\begin{proposition}\label{p5.6} If $(N, g)$ is an EALE space, then 
\begin{equation}\label{e5.26}
H_{2}(N, {\mathbb F}) \neq  0,
\end{equation}
for some field $F$. 
\end{proposition}

{\bf Proof:} From the results discussed above, one has $b_{1}(N) = b_{3}(N) 
= 0$, and so 
$$\chi (N) = 1 + b_{2}(N) \geq  1. $$
If $b_{2}(N) \neq 0$, then of course \eqref{e5.26} holds with ${\mathbb F} = 
{\mathbb R}$, so suppose $b_{2}(N) = 0$, so that $\chi(N) = 1$. The Euler characteristic 
can be computed with homology with coefficients in any field ${\mathbb F}$, so that 
$$\chi (N) = 1 - H_{1}(N, {\mathbb F}) + H_{2}(N, {\mathbb F}), $$
since again $H_{3}(N, {\mathbb F}) = H_{4}(N, {\mathbb F}) = 0$. If \eqref{e5.26} 
does not hold, then it follows that also $H_{1}(N, {\mathbb F}) = 0$, for all fields 
${\mathbb F}$. By the universal coefficient theorem, this implies that 
$H_{1}(N, {\mathbb Z}) = H_{2}(N, {\mathbb Z}) = 0$. 

 This means that $N$ is an integral homology ball, with finite 
$\pi_{1}$. A standard argument using the exact sequence of the pair $(N, \partial N)$ 
and Poincar\'e-Lefschetz duality shows that $\partial N \simeq S^{3}/\Gamma$ is 
a homology 3-sphere, and hence is the Poincar\'e homology sphere, with $\Gamma$ 
equal to the binary icosahedral group of order 120. One can now use a simple 
argument based on the $\eta$-invariant of $S^{3}/\Gamma$ to obtain a contradiction, 
cf.~\cite{An2} for further details. 
{\endproof}

  Proposition 5.6 shows that EALE spaces $N$ have non-trivial 2-dimensional 
topology. By construction, such spaces are naturally embedded in $M$, via 
\eqref{e5.18} for instance. 

\begin{proposition}\label{p5.7} For any EALE space $N \subset M$, the inclusion 
$\iota: N \hookrightarrow M$ induces an injection 
\begin{equation}\label{e5.27}
0 \rightarrow  H_{2}(N, {\mathbb F}) \rightarrow  H_{2}(M, {\mathbb F}).
\end{equation}
\end{proposition}

{\bf Proof:} The Mayer-Vietoris sequence for a thickening of the pair 
$(N, M\setminus N)$ gives
\begin{equation}\label{e5.28}
H_{2}(S^{3}/\Gamma , {\mathbb F}) \rightarrow  H_{2}(N, {\mathbb F})\oplus 
H_{2}(M\setminus N, {\mathbb F}) \rightarrow  H_{2}(M, {\mathbb F}).
\end{equation}
It suffices then to show that the inclusion map $i: \partial N = 
S^{3}/\Gamma  \rightarrow  N$ induces 0 on homology, i.e.~the map
$$H_{2}(S^{3}/\Gamma , {\mathbb F}) \rightarrow  H_{2}(N, {\mathbb F}) $$
is the zero map. If $N$ is simply connected, then this is clear since 
$H_{2}(S^{3}/\Gamma, {\mathbb F})$ is torsion while $H_{2}(N, {\mathbb F})$ 
is torsion-free, by Proposition 5.5 and the universal coefficient theorem. 
In general, suppose $\Sigma$ is an essential 2-cycle in $S^{3}/\Gamma$; one 
needs to prove it bounds a 3-chain in $N$. Any 2-cycle in $S^{3}/\Gamma$ 
with coefficients in ${\mathbb F}$ may be represented by a collection of 
maps $f: S^{2} \rightarrow S^{3}/\Gamma$. Let $\pi: \widetilde N \rightarrow 
N$ be the universal cover of $N$, so that $\pi$ is a finite cover. Also, 
let $\partial\widetilde N = S^{3}/\widetilde \Gamma$, so that $\pi$ induces 
a map $\pi_{\partial}: S^{3}/\widetilde \Gamma \rightarrow  S^{3}/\Gamma = 
\partial N$. The map $f$ lifts to $\widetilde f: S^{2} \rightarrow 
S^{3}/\widetilde \Gamma$. As noted above, $\widetilde f$ bounds a 3-chain 
in $\widetilde N$. Composing this with the projection map $\pi$ shows 
that $f$ bounds a 3-chain in $N$, as required. 
{\endproof}

 Via Proposition 5.6, one sees that orbifold limits in $\bar{\mathcal E}$ 
necessarily crush essential $2$-cycles in $N$ to points. This leads 
easily to the following:

\begin{corollary} \label{c5.8} Suppose 
\begin{equation}\label{e5.29}
H_{2}(M, {\mathbb F} ) = 0,
\end{equation}
for any field ${\mathbb F}$. Then $\bar{\mathcal E}_{GH} = {\mathcal E}$, 
i.e.~the moduli space ${\mathcal E}$ is complete and locally compact in 
the Gromov-Hausdorff, (or finite diameter), topology. 
\end{corollary}

{\bf Proof:} If $N$ is any EALE space associated to an orbifold limit in 
 $\bar{\mathcal E}_{GH}$, then Proposition 5.7 and \eqref{e5.29} imply that 
$H_{2}(N, {\mathbb F}) = 0$, for any ${\mathbb F}$, which contradicts 
\eqref{e5.25}. Hence, there are no such EALE spaces, and the result follows. 
{\endproof}

 The condition \eqref{e5.29} holds for instance for $S^{4}$, with any 
(potentially exotic) differentiable structure or more generally any 
integral homology sphere with torsion-free $\pi_{1}$; there are many 
such $4$-manifolds $M$, cf.~\S 7.

 It also follows from Propositions 5.6 and 5.7 that when \eqref{e5.29} holds, 
there are only finitely many components to the part 
$${\mathcal E}_{\lambda_{0}}\subset  {\mathcal E}$$
of the moduli space for which $\lambda \geq \lambda_{0} > 0$, and further 
each such component is compact in the $C^{\infty}$ topology. While this seems 
interesting and useful, we know of no applications this result. For 
instance, can it be useful in deciding whether $S^{4}$ has any Einstein 
metrics of positive scalar curvature besides the round metric?

 Without the homology condition \eqref{e5.29}, it is an open question whether 
$\bar{\mathcal E}_{\lambda_{0}}$ for $\lambda_{0} > 0$ is compact or not. For 
example, it may apriori be possible that there exists a sequence 
$g_{i}$ in distinct components of ${\mathcal E}_{\lambda_{0}}$ which 
converges to an Einstein orbifold metric in a limiting component of 
$\bar{\mathcal E}_{GH}$, (so that ${\mathcal E}_{\lambda_{0}}$ has infinitely 
many components). This is related to the following question.

\medskip

 {\it Question}. Let $(M, g_{i})$ be a sequence of Einstein metrics on $M$ 
converging to an Einstein orbifold metric $(V, g)$ associated to $M$. Does 
there exist a (continuous) curve $\gamma (t)$, $t\in (0,1]$, of Einstein 
metrics on $M$ such that, for $i$ sufficiently large, $\gamma (t_{i}) = 
g_{i}$, for some sequence $t_{i} \rightarrow 0$?

 If the answer to the question is yes, then it follows that 
${\mathcal E}_{\lambda_{0}}$ has finitely many components for 
$\lambda_{0} > 0$. Another approach toward the compactness of 
$\bar{\mathcal E}_{\lambda_{0}}$ is the following:

\medskip

  {\it Question}. Is the completion $\bar{\mathcal E}_{GH}$ a real-analytic 
variety?

   In other words, can one extend Theorem 5.1 to include orbifold formation? 
Note that in a real-analytic variety ${\mathcal V}$, an infinite sequence of 
points which converges to a limit point in ${\mathcal V}$ can be connected by 
a real-analytic curve in ${\mathcal V}$. 

  In this context, it is worth bearing in mind however that currently, except 
in the context of K\"ahler-Einstein metrics, one has no examples where the positive 
moduli space ${\mathcal E}_{\lambda_{0}}$, $\lambda_{0} > 0$, is larger than a 
single point.

\section{Moduli Spaces II.}
\setcounter{equation}{0}

 A more complete theory of the global behavior of the components of the 
moduli space ${\mathcal E}$ was developed in \cite{An3}. From a broad 
perspective, the overall picture has a strong resemblance with the 
moduli space of constant curvature metrics on surfaces, which we survey 
very briefly. 

 Recall that the moduli space of unit volume constant curvature 
metrics on an oriented surface $\Sigma$ has the following structure:

$\bullet$ $\Sigma = S^{2}$. Then ${\mathcal E} = pt$. 

$\bullet$ $\Sigma = T^{2}$. Then ${\mathcal E} = {\mathbb H}^{2} / 
SL(2,{\mathbb Z})$, the modular quotient of the upper half plane. 
Any divergent sequence $\{g_{i}\} \in {\mathcal E}$ collapses, in that 
the injectivity radius $inj_{g_{i}}(x)$ of $g_{i}$ at $x$ satisfies 
$$inj_{g_{i}}(x) \rightarrow 0,$$
for all $x\in\Sigma$. Any pointed Gromov-Hausdorff limit of $\{g_{i}\}$ 
is a line ${\mathbb R}$.

$\bullet$ $\Sigma = \Sigma_{g}$, $g \geq 2$. Then the Riemann moduli 
space ${\mathcal E}$ is a connected orbifold, of dimension $6g-6$. An 
element in $\partial{\mathcal E}$ is given by a finite collection of 
cusps, i.e.~complete hyperbolic metrics on a collection of punctured 
surfaces of total genus $g$. 

\medskip

 One may study the completion of ${\mathcal E}$ (or $\bar{\mathcal E}_{GH}$) 
in the pointed Gromov-Hausdorff topology, which allows the diameter of 
sequences of metrics in ${\mathcal E}$ to diverge to infinity. However, 
as in the case of surfaces, it is more useful to study the boundary 
$\partial{\mathcal E}$ with respect to other, more natural, metric 
topologies. In dimension 2, one of the most natural and well-studied metrics 
is the Weil-Petersson metric. This metric is just the restriction of 
the usual $L^{2}$ metric on the space $Met(\Sigma)$ of all metrics 
on $\Sigma$ to the moduli space. The $L^{2}$ metric on $Met(M)$ is 
given by
$$\langle h_{1}, h_{2} \rangle_{g} = \int_{M}\langle h_{1}(x), h_{2}(x) 
\rangle_{g}dV_{g},$$
where $h_{1}, h_{2} \in T_{g}Met(M)$. Thus, we consider the completion 
$\bar{\mathcal E}$ of ${\mathcal E}$ with respect to the $L^{2}$ metric, 
which provides more information than the completion with respect to the 
pointed Gromov-Hausdorff topology. 

 Before stating the main result, one needs the following definition. A 
domain $\Omega$ (i.e.~an open 4-manifold) weakly embeds in $M$, $\Omega  
\subset\subset M$, if for any compact subdomain $K \subset \Omega$, 
there is a smooth embedding $F = F_{K}: K \rightarrow M$. The same 
definition applies if $\Omega$ is an orbifold. 

\begin{theorem} \label{t6.1} \cite{An3}. The completion $\bar {\mathcal E}$ of 
${\mathcal E}$ with respect to the $L^{2}$ metric on ${\mathcal E}$ is a 
complete Hausdorff metric space, whose frontier $\partial{\mathcal E}$ 
consists of two parts: the orbifold part $\partial_{o}{\mathcal E}$ and 
the cusp part $\partial_{c}{\mathcal E}$. 

 {\rm I.} $\partial_{o}{\mathcal E}$ consists of Einstein orbifolds associated 
to $M$, of unit volume and finite diameter. The partial completion ${\mathcal E}
\cup\partial_{o}{\mathcal E}$ is locally compact and all the results of \S 5 hold 
as before. 

 {\rm II.} An element in the cusp boundary $\partial_{c}{\mathcal E}$ is given 
by a pair $(\Omega, g)$, where $\Omega$ is a non-empty maximal 
orbifold domain weakly embedded in $M$, consisting of a finite 
number of components $\Omega_{k}$ called cusps, each with a bounded 
number (possibly zero) of orbifold singularities. The metric $g$ is a 
complete Einstein metric on $\Omega$, with 
\begin{equation} \label{e6.1}
vol_{g}\Omega  = 1 ,
\end{equation}
and outside a compact set $K$, $\Omega$ carries an $F$-structure along 
which $g$ collapses with locally bounded curvature as one goes to 
infinity in $\Omega$; thus as $x \rightarrow \infty$ in $\Omega$, 
\begin{equation}\label{e6.2}
inj(x) \rightarrow 0 \ \ {\rm and} \ \  (|R|inj^{2})(x) \rightarrow 0. 
\end{equation}

 To describe the behavior of the region $M\setminus K$, let $g_{i}$ 
be a sequence in ${\mathcal E}$ with $g_{i} \rightarrow g \in 
\partial_{c}{\mathcal E}$ in the $L^{2}$ metric. Then $M \setminus K$ 
also carries an $F$-structure on the complement of a finite number of 
arbitrarily small balls. Thus, there exists a finite collection of 
points $z_{k}\in M$, with $dist_{g_{i}}(z_{k}, K) \rightarrow \infty$ 
as $i \rightarrow \infty$, and a sequence $\varepsilon_{i} \rightarrow 0$ 
such that outside $B_{z_{k}}(\varepsilon_{i})$, $M \setminus K$ has an 
$F$-structure. If $K_{i}$ is an exhaustion of $\Omega$, then $M\setminus K_{i}$ 
collapses everywhere as $i \rightarrow \infty$, and collapses with locally 
bounded curvature as in \eqref{e6.2} away from the singular points $\{z_{k}\}$. 

 Further, Case {\rm II}, i.e.~cusps, can form only on the components of 
${\mathcal E}$ for which there is a constant $\lambda_{0}$ such that
\begin{equation} \label{e6.3}
\lambda \leq  \lambda_{0} < 0. 
\end{equation}
\end{theorem}

  Note that the product in \eqref{e6.2} is scale invariant, so that \eqref{e6.2} 
shows that the metric becomes flat on the scale of the injectivity radius. 

 The convergence in Case I is also in the Gromov-Hausdorff topology, while 
that in Case II is also in the pointed Gromov-Hausdorff topology. One sees 
of course some similarities here with the structure of ${\mathcal E}$ in 
the case of surfaces. Note for instance that the compactness of 
$\bar {\mathcal E}_{GH}$ or $\bar {\mathcal E}$ on the components of 
positive scalar curvature resembles the compactness of Einstein metrics 
on $S^{2}$. 

\medskip

 As in the case of surfaces, (e.g.~the moduli of flat metrics on $T^{2}$), the 
$L^{2}$ completion $\bar{\mathcal E}$ is not compact in general. An important 
example is the case of Einstein metrics on the K3 surface, which for illustration 
is worth discussing in some detail. 

  By Theorem 3.2, any Einstein metric on K3 is Ricci-flat and K\"ahler, 
and so hyperk\"ahler. To any $g \in {\mathcal E}$ is associated a 3-dimensional 
subspace $P_{g} \subset H^{2}(K3, {\mathbb R}) \simeq {\mathbb R}^{22}$, given 
by the span of the K\"ahler forms associated to $g$. The intersection form is 
of type $(3, 19)$ on $H^{2}$, and restricts to a positive definite form on 
$P_{g}$. This assignment gives the period map
$$P: {\mathcal E} \rightarrow \Gamma / G_{3.19}^{+},$$
where $G_{3,19}^{+} \simeq SO(3,19)/(SO(3)\times SO(19))$ is the 
Grassmannian of positive 3-planes in $H^{2}$ and $\Gamma$ is the 
integral lattice $SO(3,19:{\mathbb Z})\simeq Aut (H^{2}, I)$. 
The space $\Gamma / G_{3.19}^{+}$ is a 57-dimensional orbifold and 
the local Torelli theorem for K3 surfaces implies $P$ is a local 
diffeomorphism.

  However, the period map $P$ is not surjective. Namely, let $\Delta 
\subset H^{1,1}(K3, {\mathbb R})\cap H^{2}(K3, {\mathbb Z})$ be the 
set of roots, i.e.~the class of effective divisors $d$ of 
self-intersection $-2$, $(d, d) = -2$. Then $Im P$ is contained 
in the open subset of $\Gamma / G_{3.19}^{+}$ for which, for any 
$d \in \Delta$, $\omega(d) \neq 0$, for some K\"ahler form $\omega 
\in P_{g}$. 

  One has $\bar{\mathcal E}_{K3} = {\mathcal E}_{K3}\cup\partial_{o}{\mathcal E}_{K3}$, 
and the set of Einstein orbifold metrics $\partial_{o}{\mathcal E}_{K3}$ is a 
countable collection of codimension 3 varieties in $\bar {\mathcal E}_{K3}$, 
corresponding to the locus where $\omega(d) = 0$, for some K\"ahler form 
$\omega$ and $d \in \Delta$. It is proved in \cite{An3} that the period map $P$ 
extends continuously to $\bar{\mathcal E}_{K3}$ and that the extended period 
map 
\begin{equation}\label{e6.4}
\bar P: \bar{\mathcal E}_{K3} \rightarrow \Gamma / G_{3.19}^{+},
\end{equation}
is an isometry between the $L^{2}$ metric on $\bar{\mathcal E}_{K3}$ 
and the complete, non-compact, finite volume locally symmetric metric 
on $\Gamma / G_{3.19}^{+}$. 

\medskip

 Returning to the discussion in general, the behavior of $\bar{\mathcal E}$ 
at infinity is described as follows:

\begin{theorem}\label{t6.2} \cite{An3}. Suppose $g_{i}$ is a divergent sequence in 
$\bar{\mathcal E}$ such that 
\begin{equation} \label{e6.5}
\lambda_{g_{i}} \rightarrow  0, \ \ {\rm as} \ \ i \rightarrow  \infty . 
\end{equation}
Then $\{g_{i}\}$ collapses everywhere, i.e.~$inj_{g_{i}}(x) \rightarrow 0$, 
$\forall x \in M$. The collapse is along a sequence of $F$-structures 
${\mathcal F}_{i}$ and with locally bounded curvature \eqref{e6.2} metrically 
on the complement of finitely many singular points $\{z_{k}\}$. 

 Suppose instead that $g_{i}$ is a divergent sequence in $\bar{\mathcal E}$ 
such that 
\begin{equation} \label{e6.6}
\lambda_{g_{i}} \leq  \lambda_{0} <  0, \ \ {\rm as}\ \  i \rightarrow \infty . 
\end{equation}
Then $\{g_{i}\}$ either has the same behavior as above in \eqref{e6.5}, or as in 
Case {\rm II} (cusps) of Theorem 6.1, where $\Omega$ may instead have possibly 
infinitely many components, of total volume at most 1. 
\end{theorem}

  Note that in the case of \eqref{e6.5}, the collapsing singularities $\{z_{k}\}$ 
must exist. Namely, if $\{z_{k}\} = \emptyset$, then it follows from Theorem 6.2 
that the manifold $M$ admits an $F$-structure. However, it is easy to see that 
any manifold $M$ with an $F$-structure has vanishing Euler characteristic, 
$\chi(M) = 0$, cf.~\cite{CG}. By Theorem 4.1, this can happen only if $(M, g)$ 
is flat. Thus, one expects that small neighborhoods of $\{z_{k}\}$ account for 
all of $\chi(M)$. 

 Recently, building on the work in \cite{An3}, Cheeger-Tian have 
improved Theorem 6.2, and proved the following, (answering a conjecture 
in \cite{An3}):

\begin{theorem}\label{et6.3} \cite{CT}. Suppose that $g_{i}$ is a divergent 
sequence in $\bar{\mathcal E}$ such that 
\begin{equation}\label{e6.7}
\lambda_{g_{i}} \leq  \lambda_{0} <  0, \ \ {\rm as} \ \  i \rightarrow  
\infty .
\end{equation}
Then $\{g_{i}\}$ cannot collapse everywhere, and has the same behavior as 
described in Case {\rm II} of Theorem 6.1. Moreover, the collapse in \eqref{e6.2} 
is with uniformly bounded curvature away from a finite number of singular points. 
\end{theorem}

 We briefly describe the basic idea in the proof of this non-collapse result. 
The proof is by contradiction, and so suppose $\{g_{i}\}$ is a divergent 
sequence in $\bar{\mathcal E}$ (of unit volume) which collapses everywhere, 
i.e.~$inj_{g_{i}}(x) \rightarrow 0$, for all $x$. By Theorem 6.1, 
there is a finite number of singularities $\{z_{k}\}$ outside of 
which the metrics $g_{i}$ are collapsing $M$ everywhere along a sequence 
of $F$-structures with locally bounded curvature. Let $U_{i} = 
(M\setminus \cup B_{z_{k}}(\varepsilon), g_{i})$, where $\varepsilon > 0$ 
is (very) small. As above, $\chi (U_{i}) = 0$. By the Chern-Gauss-Bonnet 
theorem for manifolds with boundary, one then has
\begin{equation}\label{e6.8}
\frac{1}{8\pi^{2}}\int_{U_{i}}|R|^{2} + \int_{\partial U_{i}}Q = \chi 
(U_{i}) = 0.
\end{equation}
Here the boundary term $Q$ consists of two terms; one of the form 
$A^{3}$, where $A^{3}$ is cubic in the eigenvalues of the $2^{\rm nd}$ 
fundamental form $A$ of $\partial U_{i} \subset  U_{i}$ and a second, 
$R(A)$, which is linear in the curvature of $g_{i}$ and $A$. Now suppose 
one can prove that there is a constant $K < \infty$ such that
\begin{equation}\label{e6.9}
|Q| \leq  K.
\end{equation}
Since one can easily arrange that $vol\partial U_{i} \rightarrow 0$, it then 
follows from \eqref{e6.8} and \eqref{e6.9} that
$$\frac{1}{8\pi^{2}}\int_{U_{i}}|R|^{2} \rightarrow  0,$$
as $i \rightarrow \infty$. However, $vol U_{i} \rightarrow 1$ while 
by the assumption \eqref{e6.7}, $|R| \geq c_{0} > 0$. This gives a 
contradiction, showing that it is not possible that $\{g_{i}\}$ 
collapses everywhere. 

   Given this idea, the main point is then to prove that the collapse away 
from $\{z_{k}\}$ is necessarily with uniformly bounded curvature, (in place 
of locally uniformly bounded curvature as in \eqref{e6.2}), from which 
\eqref{e6.9} then follows easily. 
{\endproof}

 To complete the analogy with the case of surfaces, it is natural to 
conjecture, (cf.~\cite{An3}), that there are {\it no} divergent sequences in 
$\bar{\mathcal E}$ when \eqref{e6.7} holds, i.e.~\eqref{e6.7} implies 
that $\bar{\mathcal E}$ is compact. This remains currently an open problem. 
We recall here that the completion of the moduli space of hyperbolic metrics 
on a surface with respect to the Weil-Petersson metric is compact, and 
agrees with the Deligne-Mumford compactification. 

\smallskip

  The results described above lead of course to many new questions, and 
we discuss some of these next. 

   First, while the process leading to the formation of orbifold singularities 
in $\bar{\mathcal E}_{GH}$, as well as the structure of the Einstein orbifold 
metrics themselves, is comparatively well understood, almost nothing is known 
about the structure of the singular points $\{z_{k}\}$ arising in collapsing 
regions. One would expect that such singularities can be modelled on spaces 
such as the multi-Taub-NUT metrics \eqref{e5.22} and other ALF spaces, as in 
\cite{CK} for instance; however, with the exception of the work of Gross-Wilson 
\cite{GW} discussed in \S 7, there are no studies along these lines. 

   For instance, the fact that the $F$-structure does not extend through 
the singularities $\{z_{k}\}$ suggests the following question, which is 
of independent interest:

\smallskip

 {\it Question}. Suppose $(N, g)$ is a complete, non-compact Ricci-flat 
4-manifold, with a free isometric $S^{1}$ action. Is $(N, g)$ necessarily 
flat?

 This is of course false if the $S^{1}$ action is allowed to have a 
non-trivial fixed point set; for example, the Gibbons-Hawking metrics 
\eqref{e5.20} have an isometric $S^{1}$ action. This result is also easy to 
prove if $(N, g)$ is ALE for instance. Moreover, there is a positive answer to 
the question in the setting of Ricci-flat Lorentz metrics with a non-vanishing 
time-like Killing field, cf.~\cite{An5}. However, the method of proof used there 
does not carry over to the Riemannian case.  

   Similarly, very little is known in detail about the possible formation of 
Einstein cusp manifolds as described in Theorem 6.1, beyond the obvious 
examples of this behavior on products of surfaces of genus $> 2$, (and 
the K\"ahler-Einstein case on surfaces of general type); see however \S 7 
for some further results. 

   For instance, it is not known if the cusp ends must be of finite 
topological type. A typical cross-section of a cusp end is a closed 
3-manifold which collapses with bounded curvature. These are graph 
manifolds, i.e.~unions of Seifert fibered 3-manifolds along tori. 
Can general graph manifolds arise in the cusp ends? One has the 
standard examples of 3-tori $T^{3}$ and the 3-dimensional Nil manifolds, 
arising from cusps of hyperbolic and complex hyperbolic space-forms, as 
well as collapse along circles in the case of products of hyperbolic 
surfaces; however, beyond this, no further examples are known. 
 
  As discussed in \S 5, orbifold singularities can only arise when essential 
2-cycles in $M$ are crushed to points, cf.~Propositions 5.5-5.8. No analog of 
this result is known for cusps, and it would be very interesting to show 
cusps or their complements are topologically essential in $M$ in some way. 

  Another very interesting question arises again by comparison with the case 
of surfaces. For surfaces, the Deligne-Mumford or Weil-Petersson compactification 
of the Riemann moduli space $\bar{\mathcal E}$ is a {\it closed} complex analytic 
variety, for which the boundary $\partial {\mathcal E} = \bar{\mathcal E} \setminus 
{\mathcal E}$ is a complex subvariety of real codimension 2. Thus, the boundary 
is not ``topological'', in the sense that $\bar{\mathcal E}$ is a (singular) 
manifold with boundary, (with boundary of codimension 1). Is the same true for 
the $L^{2}$ compactification $\bar{\mathcal E}$ in dimension 4, for the 
components with $\lambda < 0$? In other words, can one give the completion 
$\bar{\mathcal E}$ the structure of a cycle?

\section{Constructions of Einstein metrics, II.}
\setcounter{equation}{0}

  In this section, we discuss a different approach to constructing 
Einstein metrics, namely by ``glueing'' or resolving singular Einstein 
metrics. Beginning with the work of Taubes \cite{Ta}, this singular perturbation 
method has been very successful in constructing interesting solutions of a wide 
variety of geometric equations and it is of interest to study this 
issue in the context of Einstein metrics. 

   We have seen in prior sections that one can describe in some detail 
the singular Einstein structures that form at the boundary or at infinity 
of the moduli space ${\mathcal E}$. The singular behaviors are of three types.

$\bullet$ orbifold singularities.

$\bullet$ collapsing singularities or limits.

$\bullet$ cusp limits.

  The collapsing singularities are point-type singularities whose 
blow-ups are modeled on complete Ricci-flat spaces which are not ALE; 
for example the ALF Taub-NUT or multi-Taub-NUT metrics discussed in 
\S 5. The collapsed or cusp limits describe global limit behavior, as 
opposed to local singularities.

   For each of these types of singular behavior, the natural glueing question 
is: when can this process be reversed, i.e.~given in general such a singular Einstein 
metric or configuration, when can one ``resolve the singularity'' by finding 
a sequence or curve of Einstein metrics, on a given compact manifold $M$ or 
sequence of compact manifolds $M_{i}$, which tend to the given singular space or 
configuration in the limit.

   As with all glueing procedures, the construction proceeds in two steps, 
one more conceptual and one more technical. First, one constructs an 
(essentially explicit) good approximate solution to the Einstein equations 
on $M$, (or $M_{i}$), by patching together exact Einstein metrics on the 
domains which together cover $M$. Next one needs to prove that such an 
approximate solution can be perturbed to a nearby exact solution of the 
Einstein equations. This perturbation to an exact solution is carried out 
by means of the inverse or implicit function theorem, and often requires 
a considerable amount of technical work to establish. 

\smallskip

{\bf I. Orbifold singularities.} 

  It has long been an open question whether Einstein orbifold metrics 
$(V, g_{\infty})$ can be resolved to smooth Einstein metrics $(M, g)$ 
close to $(V, g_{\infty})$ in the Gromov-Hausdorff topology. The idea 
here would be to reverse the process of formation of orbifold singularities 
described in \S 5. Here there is often no difficulty in carrying out 
the first step in the construction. For instance, if a singularity 
$q$ of $V$ is of the type $C(S^{3}/\Gamma)$ with $\Gamma \subset SU(2)$, 
then it is easy to construct good approximate solutions by gluing the 
punctured Einstein orbifold $V \setminus \{q\}$ with a truncation and 
rescaling of a hyperk\"ahler EALE space given by Theorem 5.4. The main 
issue is then the perturbation to an exact Einstein metric. 

   There are certainly some situations where one can resolve orbifold 
singularities. Notably, this is the case on K3 surfaces. As a concrete 
example, cf.~\cite{Pa2}, consider the classical Kummer construction of K3 
surfaces. Thus, let $T^{4} = {\mathbb C}^{2}/\Lambda$ be a complex torus, 
and consider the involution $A: T^{4} \rightarrow T^{4}$, $A(x) = -x$. 
This map has 16 fixed points, (the points of order 2), and the quotient 
$V = T^{4}/{\mathbb Z}_{2}$ is an orbifold with 16 singular points, 
each of the form $C({\mathbb R}{\mathbb P}^{3})$. The flat metric 
$g_{0}$ is an orbifold singular Einstein metric on $V$. One may then 
take 16 copies of the Eguchi-Hanson metric \eqref{e5.21}, truncated and 
scaled down, and glue this onto the regular set $(V_{0}, g_{0})$. This gives 
an approximate Ricci-flat metric on K3. 

  It follows from the discussion on the K3 surface in \S 6 that there 
are smooth Einstein metrics on the K3 surface very close to such 
appoximate solutions, cf.~also \cite{KT}. In fact, by \eqref{e6.4}, all the 
orbifold Einstein metrics associated with K3 can be resolved. However, 
this proof is indirect, relying on the structure of $\bar{\mathcal E}_{K3}$, 
and does not lead to the actual construction of new Einstein metrics. 

  Moreover, the resolution of orbifold singularities in this way 
does not work in general for K\"ahler-Einstein orbifolds. This 
follows from fact that there are many more orbifold singular 
K\"ahler-Einstein metrics than smooth K\"ahler-Einstein metrics, 
cf.~\cite{Ko1,Ko2} for instance. For example, consider log del Pezzo 
surfaces; these are surfaces with $c_{1} > 0$ with quotient singularities, 
i.e.~singularities of the form ${\mathbb C}^{2}/\Gamma$, $\Gamma 
\subset U(2,{\mathbb C})$. It is proved in \cite{Ko1} that there are log 
del Pezzo surfaces $V$ having orbifold K\"ahler-Einstein metrics, with $\lambda 
> 0$, which have arbitrarily large $b_{2}$. However, although all smooth del 
Pezzo surfaces have K\"ahler-Einstein metrics with $\lambda > 0$, 
(by Theorem 3.3), all such satisfy $b_{2} \leq 9$. 

  If one could carry out a glueing construction in the K\"ahler-Einstein context, 
resolving the singularities of $V$, then a smooth K\"ahler-Einstein metric would 
also have $\lambda > 0$, and hence is defined on a smooth del Pezzo surface $M$, 
as in Theorem 3.3. The same argument establishing \eqref{e5.28} also proves 
$H_{2}(M\setminus N, {\mathbb R})$ injects in $H_{2}(M, {\mathbb R})$, which 
gives a contradiction. 

  While it is still conceivable that such orbifold metrics could be resolved 
by (real) Einstein metrics, this seems unlikely. It would of course be very 
interesting to establish this one way or the other.  

\smallskip

{\bf II. Collapsing configurations.} 

   An interesting glueing result was proved by Gross-Wilson \cite{GW}, describing 
the collapse behavior of (some) sequences of Einstein metrics tending to 
infinity in the moduli space ${\mathcal E}_{K3}$ of Einstein metrics on 
the K3 surface. 

   To describe the result, a large family of K3 surfaces are given as 
elliptic (i.e.~torus) fibrations; there is a holomorphic map $K3 \rightarrow S^{2}$, 
with 24 singular fibers, each of type $I_{1}$, (a pinched torus). The approximate 
solution is given by glueing together two types of Ricci-flat K\"ahler metrics. 
First, away from the singular fibers, one takes a ``semi-flat'' metric, defined 
by a Riemannian submersion to a base given by a disc $D^{2}\subset S^{2}$ with 
flat metric induced on the fibers. The fibers are scaled to be highly collapsed, 
i.e.~of small diameter. In a neigbhorhood of each singular fibers, one takes a 
suitable truncation and rescaling of the Ooguri-Vafa metric \cite{OV}, 
i.e.~a periodic Taub-NUT metric. Thus, choose the potential $V$ as in 
\eqref{e5.22}, with infinitely many points $p_{i}$ of distance $\varepsilon$ 
apart along the $z$-axis in ${\mathbb R}^{3}$. By adding a suitable constant, 
the potential may be renormalized to sum to a smooth, periodic harmonic function. 
The resulting metric as in \eqref{e5.20} is then also periodic, and taking 
the ${\mathbb Z}$ quotient gives a metric on a manifold topologically 
equivalent to a neighborhood of a singular fiber, and which, after rescaling, 
is close to a semi-flat metric. One then perturbs this approximate solution 
to an exact Ricci-flat K\"ahler metric by obtaining uniform estimates for 
the complex Monge-Ampere equation \eqref{e3.3}. 

  Although this glueing result does not actually construct any new 
Einstein metrics, (they are already given by Yau's Theorem 3.2 on the 
K3 surface), the construction is important in understanding various 
aspects of mirror symmetry on K3. 

\medskip

  Next, we describe a much simpler glueing method, which, if successful, 
would lead to many interesting, new Einstein metrics. We will describe the 
method only in a special case, since this illustrates the main idea. 
Thus, let $(M, g)$ be a Ricci-flat manifold of the form $M = S^{1}\times N$, 
with $g$ a product metric. Assume that the metric $g$ is highly 
collapsed, so that the length of $S^{1}$ is small. Let $\hat M$ be the 
manifold obtained by performing surgery to kill the $S^{1}$ factor in 
$\pi_{1}(M)$. Thus, one removes $S^{1}\times B^{n-1} \subset M$ 
and glues in $D^{2}\times S^{n-2}$. 

  In \cite{An4}, a good approximate Einstein metric was constructed on 
$\hat M$ in that, for any given $\varepsilon > 0$, there are metrics 
$g_{\varepsilon}$ on $\hat M$ such that  
\begin{equation}\label{e7.1}
|Ric_{g_{\varepsilon}}| \leq \varepsilon.
\end{equation}
To see this, consider first the family of Schwarzschild metrics $g_{Sch}$:
\begin{equation}\label{e7.2}
g_{Sch} = V^{-1}dr^{2} + Vd\theta^{2} + r^{2}g_{S^{n-2}},
\end{equation}
where $V(r) = 1 - \frac{2m}{r^{n-3}}$ and $\theta \in [0,\beta]$; here 
$\beta = \beta(m)$ is chosen so that the metric is smooth at the horizon 
where $V = 0$. This metric lives on the manifold $D^{2}\times S^{n-2}$, and 
is asymptotic to the product metric on $S^{1}\times {\mathbb R}^{n-1}$. 
Truncating and rescaling the metric to a small size, one may then glue 
this into $\hat M \setminus (S^{1}\times B^{n-1})$ and construct 
$g_{\varepsilon}$ satisfying \eqref{e7.1}. 

 {\it Question}. Are there situations where the approximate solutions 
$g_{\varepsilon}$ on $\hat M$ can be perturbed to a nearby exact Einstein 
metric. 
 
  If so, one would expect to be able to repeat this surgery arbitrarily many 
times. 

\medskip

{\bf III. Cusp configurations.}

  Here we describe a construction of Einstein metrics obtained by ``resolving'' 
cusp configurations. The results are closely analogous to the Thurston theory 
\cite{Th1} of Dehn surgery on hyperbolic 3-manifolds. The construction holds in all 
dimensions $n \geq 3$ but for simplicity we work with $n = 4$. 

  One starts with any complete, non-compact hyperbolic 4-manifold 
$N = N^{4}$ of finite volume, with metric $g_{-1}$. The manifold $N$ has 
a finite number of cusp ends $\{E_{j}\}$, $1 \leq j \leq k$, each 
diffeomorphic to ${\mathbb R}\times F$, where $F$ is a compact flat 3-manifold. 
By passing to covering spaces if necessary, it may and will be assumed that each 
$F$ is a 3-torus $T^{3}$. 

  Now perform Dehn filling on each of the cusp ends. Thus, fix a flat torus $T^{3}$ 
in a given end $E$ and let $\sigma$ be any simple closed geodesic in $T^{3}$. Attach 
a 4-dimensional solid torus $D^{2}\times T^{2}$ onto $T^{3}$ by a diffeomorphism 
$\partial (D^{2}\times T^{2}) = T^{3}$, sending the $S^{1} = \partial D^{2}$ 
onto $\sigma$. Carrying out this Dehn filling for each of the cusp ends of $N$ gives 
a compact manifold $M_{\sigma}$, $\sigma  = (\sigma_{1}, \cdots , \sigma_{k})$. We 
will say that $\sigma$ is sufficiently large if the length $l(\sigma_{j})$ of each 
closed geodesic $\sigma_{j} \subset T_{j}^{3}$ is sufficiently long. 

\begin{theorem}\label{t7.1} \cite{An6}. For any choice of $\sigma$ sufficiently 
large, the manifold $M_{\sigma}$ admits an Einstein metric with negative scalar 
curvature. 
\end{theorem}

 To give an idea of how many Einstein metrics are constructed in this 
way, it is known from results of \cite{Bu} that the number ${\mathcal H}(V)$ 
of complete non-compact hyperbolic 4-manifolds of finite volume with 
volume $\leq  V$ grows super-exponentially: in fact
$$e^{aV\ln V}\leq  {\mathcal H} (V) \leq  e^{bV\ln V}. $$
With each such $N$, Theorem 7.1 associates infinitely many 
homeomorphism types of compact manifolds $M_{\sigma}$. Formally, the 
number of such compact manifolds is $\infty^{q}$, where $q$ is the 
number of cusp ends. (The number of cusp ends also grows linearly 
with $V$). Most of these Einstein metrics are not locally isometric 
(while all hyperbolic manifolds are locally isometric). 

  The Euler characteristic and signature of $M_{\sigma}$ are given by 
$$\chi (M_{\sigma}) = \chi (N), \ \ \tau (M_{\sigma}) = 0. $$
Each $M_{\sigma}$ is aspherical, i.e.~a $K(\pi, 1)$, (for $\sigma$ sufficiently 
large), and in fact admits metrics of non-positive, (but not negative), sectional 
curvature. A surprising result of \cite{RT} shows that there exists $N$ and infinitely 
many choices of $\sigma$ for which $M_{\sigma}$ is an integral homology 4-sphere. 

\medskip

  It is worth describing in some detail how the approximate Einstein metric is 
constructed on $M_{\sigma}$. First, one has the given hyperbolic, and so Einstein, 
metric on $N$. One needs then an Einstein metric on the solid torus $D^{2}\times T^{2}$ 
which closely matches the hyperbolic metric on a cusp end of $N$. A model 
for such metrics was constructed long ago by physicists: this is the family 
of toral AdS black hole metrics
$$g_{BH} = V^{-1}dr^{2} + Vd\theta^{2} + r^{2}g_{T^{2}}, $$
where $g_{T^{2}}$ is any flat metric on the torus and $V = V(r)$ is 
given by
$$V = r^{2} - \frac{2m}{r}, $$
compare with the Schwarzschild metric in \eqref{e7.2}. The set $r = r_{+} = 
(2m)^{1/3}$ where the potential $V$ vanishes is called the horizon $H$. 
Note that $H$ is a totally geodesic and flat torus $T^{2}$ in $g_{BH}$. 
For the metric $g_{BH}$ to be smooth at $H$, one requires that the circular 
parameter $\theta$ runs over the interval $\theta\in [0,\beta]$, where 
$\beta  = 4\pi /3r_{+}$. 

 This metric is asymptotically hyperbolic, in that the curvature tends 
to -1 at infinity, but the metric has infinite volume, and so is not at all 
close to a hyperbolic cusp metric on $N$. However, one can take quotients of 
$g_{BH}$ to obtain metrics which are almost cusp-like in large regions. 

 To see this more clearly, consider first the universal cover $D^{2}\times 
{\mathbb R}^{2}$ of the solid torus $D^{2}\times T^{2}$ with lifted metric
$$\widetilde g_{BH} = V^{-1}dr^{2} + Vd\theta^{2} + r^{2}(dx_{1}^{2} + 
dx_{2}^{2}). $$
The metrics $g_{BH}$ depend on the mass parameter $m$, but all metrics 
$\widetilde g_{BH}$ are isometric, as one easily seems by the change of 
variable $r \rightarrow  r_{m} = mr$; thus for convenience, set $m = 
\frac{1}{2}$, so that $r_{+} = 1$. Let $D(R) = \{r \leq  R\} \subset (D^{2}\times 
{\mathbb R}^{2}, \widetilde g_{BH})$ and let $S(R) = \partial D(R) 
= \{r = R\}$. The induced metric on the boundary $S(R)$ is then a flat metric
$$V(R)d\theta^{2} + (dz_{1}^{2} + dz_{2}^{2}), $$
on $S^{1}\times {\mathbb R}^{2}$, where $z_{i} = Rx_{i}$. Choose $R$ so that
$$\sqrt{V(R)}\beta  = l(\sigma), $$
so that the length of $S^{1}\times \{pt\} \subset S(R)$ equals $l(\sigma)$. 
Now given the flat stucture $g_{0}$ on $T^{3} \subset N$, observe that 
up to conjugacy there is a unique free isometric ${\mathbb Z}^{2}$ action 
on the flat product $S(R) = S^{1}\times {\mathbb R}^{2}$ such that the 
projection map to the orbit space
$$\pi : S^{1}\times {\mathbb R}^{2} \rightarrow  T^{3} $$
satisfies $\pi (S^{1}) = \sigma$, and for which the flat structure on 
$T^{3}$ induced by $\pi$ is the given structure $g_{0}$. In fact, the 
map $\pi$ is just the covering space of $(T^{3}, g_{0})$ corresponding 
to the subgroup $\langle \sigma \rangle \subset  \pi_{1}(T^{3})$. 

 One may easily verify that this ${\mathbb Z}^{2}$ action extends radially 
to a smooth and free isometric action on the domain $D(R) \subset  
D^{2}\times {\mathbb R}^{2}$. The quotient space $(D^{2}\times 
{\mathbb R}^{2})/{\mathbb Z}^{2} \simeq  D^{2}\times T^{2}$ gives 
the twisted toral AdS black hole metric
$$\hat g_{BH} = [V^{-1}dr^{2} + Vd\theta^{2} + 
r^{2}g_{{\mathbb R}^{2}}]/{\mathbb Z}^{2}. $$
Note that diameter of the core $T^{2}$ at $r = r_{+}$ has size about 
$R^{-1}$, and so is very small. The locus $\{r = R\}$ with the metric 
induced from $g_{BH}$ is now isometric to $(T^{3}, g_{0})$. Since this 
is also the metric at the boundary of the torus-truncated hyperbolic 
manifold, the boundaries may be identified, giving a $C^{0}$ metric on 
$M_{\sigma}$. The $2^{\rm nd}$ fundamental forms also are very close 
(for $R$ large), as are the curvatures. Smoothing the seam at the 
identification locus gives the approximate solution of the Einstein 
equations. 

   One then shows via the inverse function theorem that the 
approximate metric may be perturbed to a nearby exact Einstein metric 
on $M_{\sigma}$. As a consequence, it follows that the constructed 
Einstein metric $g_{\sigma}$ on $M_{\sigma}$ is locally rigid, 
i.e.~is an isolated point in the full moduli space ${\mathcal E}$ 
on $M_{\sigma}$. This local rigidity is essentially a feature of negative 
curvature. Note that it could not hold in the context of the conjectured 
Ricci-flat Schwarzschild glueings in \eqref{e7.1}-\eqref{e7.2}, which would 
give rise to non-trivial curves in the moduli space. It is unknown if 
$g_{\sigma}$ is the unique Einstein metric on $M_{\sigma}$.  

  The homeomorphism type of $M_{\sigma}$ is essentially uniquely determined 
by the choice of $[\sigma_{j}] \in \pi_{1}(T_{j}^{3})$, $1 \leq j \leq k$. 
The Einstein metrics $(M_{\sigma}, g_{\sigma})$ are all close to the 
original hyperbolic manifold $(N, g_{-1})$ in the pointed 
Gromov-Hausdorff topology, and are smoothly close away from the horizon 
region. On any sequence where $l(\sigma_{j}) \rightarrow \infty$ for all 
$j$, one has $(M_{\sigma}, g_{\sigma}) \rightarrow (N, g_{-1})$ in the 
pointed Gromov-Hausdorff topology, i.e.~such infinite sequences tend toward 
the original ``singular'' cusp configuration. Note that the original tori 
$T_{j}^{3} \subset N$, which are also embedded in each $M_{\sigma}$, are now 
no longer incompressible in $M_{\sigma}$; only the core tori $T_{j}^{2}$ are 
topologically essential in $M_{\sigma}$. 

  It is well-known that $\chi(N)$ can take on arbitrary values in 
${\mathbb Z}^{+}$, although $\chi(N)$ is even if and only if $N$ is 
orientable. Recall that for hyperbolic manifolds of finite volume, one has the 
relation, (via \eqref{e4.2}),
$$vol N = \frac{4\pi^{2}}{3}\chi (N). $$
It then follows again from \eqref{e4.2} that if the Einstein metric $g_{\sigma}$ 
on $M_{\sigma}$ is scaled so that $\lambda = -3$, then
$$vol M_{\sigma} <  vol N,$$
for any $\sigma$, (large). Thus, all manifolds obtained by performing Dehn 
filling on the ends of $N$ have volume strictly less than $vol N$, as in 
the Thurston theory of Dehn surgery. 

  This cusp formation on limits of sequences does not take place within 
the moduli space ${\mathcal E}$ on a fixed manifold. In fact, with the 
exception of the case K\"ahler-Einstein metrics with $c_{1} < 0$, there 
are no examples of {\it curves} of Einstein metrics in a fixed component 
of ${\mathcal E}$ which limit on a cusp configuration. It would be very 
interesting to find examples of such curves. Similarly, it would be 
interesting to understand if such ``cusp resolution'' can be carried 
out on other types of cusps, for instance on complex hyperbolic cusps, 
or cusps arising from products of surfaces. 

  It is also worth pointing out that this process of Dehn surgery is an 
important method of constructing exotic smooth structures on a 4-manifold 
of a given homeomorphism type, cf.~\cite{FS}. It remains completely open whether 
the construction above can be carried out in this context, to give a 
construction of Einstein metrics on such exotic smooth structures.

\section{Concluding Remarks}
\setcounter{equation}{0}

  The survey above shows that we are far from any theory describing the structure 
of Einstein metrics on 4-manifolds. One has instead an interesting collection of 
different methods, ideas and results which, at the moment, do not assemble to give 
any coherent picture or structure. 

   A basic issue is to what extent the Thurston picture of 3-manifolds, proved 
by Perelman, carries over to dimension 4; thus to what extent does a general 4-manifold 
decompose into a collection of domains, each of which carries a complete Einstein, 
or ``Einstein-like'' metric, or collapses along an $F$-structure. From the natural 
viewpoint of the Ricci flow, one should allow Ricci solitons as natural generalizations 
of Einstein metrics. As a first step for instance, it would be very interesting to know 
if there is an analog of Hamilton's structure theorem \cite{Ha2} for non-singular 
solutions of the Ricci flow in dimension 4.

\bibliographystyle{plain}

\begin{thebibliography}{WWW}
\bibitem{Pe1} G. Perelman, The entropy formula for the Ricci flow and its geometric 
applications, math.DG/0211159.
\end{thebibliography}

\bigskip
\noindent
\address{Department of Mathematics\\
S.U.N.Y. at Stony Brook\\
Stony Brook, N.Y. 11794-3651}

\noindent
E-mail: anderson@math.sunysb.edu

\end{document}